\documentclass[twoside, onecolumn, 10pt]{article}
\usepackage{amssymb}
\usepackage{mathrsfs}
\usepackage{amsmath}
\newtheorem{theorem}{Theorem}[section]
\newtheorem{lemma}[theorem]{Lemma}
\newtheorem{corollary}[theorem]{Corollary}

\newtheorem{problem}[theorem]{Problem}

\def\bc{\begin{center}}
\def\ec{\end{center}}

\newcommand{\nb}{\nonumber}
\newcommand{\R}{{\mathscr R}}

\begin{document}

\noindent
{\Large \bf Equalities and inequalities for Hermitian  solutions and Hermitian definite solutions
of the two matrix equations $AX = B$ and $AXA^* = B$}\\

\noindent
YONGGE  TIAN\\

\noindent {\bf Abstract.} \ This paper studies algebraic properties of Hermitian solutions and Hermitian definite solutions of
the two types of matrix equation $AX = B$ and $AXA^* = B$.  We first establish a variety of rank and inertia formulas for calculating the maximal
and minimal ranks and inertias of Hermitian solutions and Hermitian definite solutions of the matrix
equations $AX = B$ and $AXA^* = B$, and then use them to characterize many qualities and inequalities for Hermitian
solutions and Hermitian definite solutions of the two matrix equations and their variations. \\


\medskip


\renewcommand{\thefootnote}{\fnsymbol{footnote}}
\footnotetext{\hspace*{-5mm}This work was supported partially by National Natural Science Foundation of
China (Grant No. 11271384).}

\section{Introduction}

\renewcommand{\theequation}{\thesection.\arabic{equation}}
\setcounter{section}{1}
\setcounter{equation}{0}

Consider the following two well-known
linear matrix equations
\begin{equation}
AX = B
\label{11}
\end{equation}
and
\begin{equation}
AXA^* = B,
\label{12}
\end{equation}
both of which are simplest cases of various types of linear
matrix equation (with symmetric patterns), and are the starting point of many
advanced study on complicated matrix equations. A huge amount of results about
the two equations and applications were given in the literature. In particular, many
problems on algebraic properties of solutions of the two matrix equations were explicitly
characterized by using formulas for ranks and inertias of matrices.
In this paper, the author focuses on Hermitian solutions or Hermitian definite solutions of (\ref{11}) and   (\ref{12}),
 and studies the following optimization problems on the ranks and inertias
of Hermitian solutions and Hermitian definite solutions of (\ref{11}) and  (\ref{12}):

\begin{problem} \label{T11}
{\rm
Let  $A, \, B \in {\mathbb C}^{m\times n}$ be given, and
assume that (\ref{11}) has a Hermitian solution or Hermitian definite solution $X\in {\mathbb C}^{n\times n}$.
In this case, establish formulas for calculating the extremal ranks and inertias of
\begin{equation}
X  - P,
\label{13}
\end{equation}
where $P \in {\mathbb C}^{n\times n}$ is a given Hermitian matrix, and use the formulas to characterize behaviors
of these Hermitian solution and definite solution, in particular, to give necessary and sufficient
 conditions for the four inequalities
\begin{equation}
X \succ  P, \ \  X \succcurlyeq P, \ \ X \prec P, \ \ X \preccurlyeq P
\label{14}
\end{equation}
in the L\"owner partial ordering to hold, respectively.
}
 \end{problem}

\begin{problem} \label{T12}
 {\rm
Let $A, \, B, \,  C, \, D \in {\mathbb C}^{m\times n}$
be given,  $X, \, Y \in {\mathbb C}^{n\times n}$ be two unknown matrices, and assume that the two linear matrix
equations
 \begin{equation}
AX = B,  \ \ CY = D
\label{15}
\end{equation}
 have Hermitian solutions, respectively. In this case,  establish formulas for calculating the extremal
 ranks and inertias of the difference
 \begin{equation}
X - Y
\label{16}
\end{equation}
of the Hermitian solutions, and use the formulas to derive necessary and sufficient conditions for
\begin{equation}
X \succ Y, \ \  X \succcurlyeq Y, \ \  X \prec Y, \ \ X \preccurlyeq Y
\label{17}
\end{equation}
to hold in the L\"owner partial ordering,  respectively.
}
\end{problem}

\begin{problem} \label{T13}
{\rm
Let  $A\in \mathbb C^{m\times n}$ and $B \in {\mathbb C}_{{\rm H}}^{m}$ be given, and assume that
(\ref{12}) has a Hermitian solution. Pre- and post-multiplying a
matrix $T \in \mathbb C^{p\times n}$ and its conjugate transpose $T^*$ on both sides
of (\ref{12}) yields a transformed equation as follows
 \begin{equation}
 TAXA^*T^{*} = TBT^*.
 \label{18}
\end{equation}
Further define
\begin{align}
& {\cal S} = \{\, X \in {\mathbb C}_{{\rm H}}^{n} \ | \ AXA^* = B   \, \},
\label{19}
\\
& {\cal T} = \{\,  Y \in {\mathbb C}_{{\rm H}}^{n} \ | \ TAYA^*T^* = TBT^*  \, \}.
\label{110}
\end{align}
In this case, give necessary and sufficient conditions for  ${\cal S} = {\cal T}$ to hold,
as well as necessary and sufficient conditions for
$X \succ Y,$ $X \succcurlyeq Y,$ $X \prec Y$ and $X \preccurlyeq Y$ to hold for $X\in {\cal S}$ and $Y \in  {\cal T}$,
 respectively.
}
\end{problem}

\begin{problem} \label{T14}
{\rm
Assume that  {\rm (\ref{12})} has a Hermitian solution,  and define
\begin{align}
&  {\cal S} = \left\{ \, X  \in {\mathbb C}_{{\rm H}}^{n} \, | \, AXA^* = B  \right\},
\label{111}
 \\
& {\cal T} = \left\{ \,(X_1 + X_2)/2 \, | \, T_1AX_1A^*T^*_1 = T_1BT^*_1, \, T_2AX_2A^*T^*_2 = T_2BT^*_2, \
 X_1, \, X_2 \in {\mathbb C}_{{\rm H}}^{n}  \, \right\}.
\label{112}
\end{align}
In this case, give necessary and sufficient conditions for  ${\cal S} = {\cal T}$ to hold.
}
\end{problem}

\begin{problem} \label{T15}
{\rm
Denote the sets of all least-squares solutions and least-rank Hermitian solutions of (\ref{12}) as
\begin{align}
{\cal S} & =  \left\{ \,  X \in \mathbb C^{n}_{{\rm H}} \ | \ \|\, B - AXA^{*} \,\|_F  = \min \, \right\},
\label{113}
\\
{\cal T} & =  \left\{ \, Y \in \mathbb C^{n}_{{\rm H}} \ | \  r(\, B - AYA^{*} \,)  = \min  \, \right\}.
\label{114}
\end{align}
In this case, establish necessary and sufficient conditions for $X \succ Y,$
$X \succcurlyeq Y,$ $X \prec Y$ and $X \preccurlyeq Y$ to hold for $X\in {\cal S}$ and $Y \in  {\cal T},$
 respectively.}
\end{problem}

Matrix equations have been a prominent concerns in matrix theory and applications. As is known to all,
two key tasks in solving a matrix equation is to give identifying condition for the existence of a solution of the equation,
and to give general solution of the equation. Once general solution is given, the subsequent work is
to describe behaviors of solutions of the matrix equation, such as, the uniqueness of solutions;
the norms of solutions; the ranks and ranges of solutions, the definiteness of solutions,
equalities and inequalities of solutions, etc. Problems \ref{T11} and \ref{T12} describe the inequalities for solutions
of (\ref{11}), as well as relations between solutions of two linear matrix equations.

Throughout this paper, ${\mathbb C}^{m\times n}$ and ${\mathbb C}_{{\rm H}}^{m}$
stand for the sets of all $m\times n$ complex matrices  and all $m\times m$ complex Hermitian
 matrices, respectively; the symbols $A^{*}$, $r(A)$ and ${\mathscr R}(A)$ stand for the transpose,
conjugate transpose, rank and range (column space) of a matrix $A\in
{\mathbb C}^{m\times n}$, respectively; $I_m$ denotes the identity
matrix of order $m$; $[\, A, \, B \,]$ denotes a row block matrix consisting of $A$ and $B$.
 The Moore--Penrose inverse of a matrix $A \in {\mathbb C}^{m \times n}$, denoted by $A^{\dag}$,
 is defined to be the unique matrix $X \in {\mathbb C}^{n \times m}$ satisfying the
 following four matrix equations
$$
{\rm (i)} \  AXA = A,  \ \ {\rm (ii)} \ XA X = X, \ \ {\rm (iii)} \
(AX)^{*} = AX, \ \  {\rm (iv)} \ (XA)^{*} =XA.
$$
Further, let $E_A = I_m - AA^{\dag}$ and $F_A = I_n - A^{\dag}A$, which ranks are given by
$r(E_A) = m - r(A)$ and $r(F_A) = n - r(A)$.  A well-known property of the Moore--Penrose inverse is $(A^{\dag})^{*} =
(A^{*})^{\dag}$. In particular$,$ both $(A^{\dag})^{*} = A^{\dag}$ and $AA^{\dag} =A^{\dag}A$ hold if $A$ is Hermitian,
i.e., $A = A^{*}$. $A \succcurlyeq 0$ ($A \succ  0$) means that $A$ is Hermitian positive semi-definite (Hermitian positive definite).
Two $A,  \, B \in {\mathbb C}_{{\rm H}}^{m}$ are said to satisfy the inequality $A \succcurlyeq B$
($A \succ  B)$ in the L\"owner partial ordering if $A - B$ is Hermitian positive semi-definite (Hermitian positive definite).
$i_{\pm}(A)$ denotes the numbers of the positive and negative eigenvalues of a Hermitian matrix
$A$ counted with multiplicities, respectively.

The results on ranks and inertias of matrices in Lemmas \ref{T16} and \ref{T17} below are obvious or well-known
(see also \cite{T-laa10,T-laa11} for their references), while the closed-form formulas for matrix ranks and inertias
in Lemmas \ref{T19} and \ref{T112}--\ref{T115}  were established by the present author, which we
shall use in the latter part of this paper to derive analytical solutions to Problems \ref{T11}--\ref{T15}.

\begin{lemma} \label{T16}
Let ${\mathcal S},$  ${\cal S}_1$ and ${\cal S}_2$ be three sets consisting of $($square$)$ matrices
over  ${\mathbb C}^{m\times n},$ and let ${\mathcal H}$ be a set consisting of Hermitian matrices over
${\mathbb C}_{{\rm H}}^{m}.$  Then$,$ the following hold$.$
\begin{enumerate}
\item[{\rm (a)}] Under $m =n,$ ${\mathcal S}$ has a nonsingular matrix if and only if
$\max_{X\in {\mathcal S}} r(X) = m.$

\item[{\rm (b)}] Under $m =n,$ all $X\in {\mathcal S}$ are nonsingular if and only if
 $\min_{X\in {\mathcal S}} r(X) = m.$

\item[{\rm (c)}] $0\in {\mathcal S}$ if and only if
$\min_{X\in {\mathcal S}} r(X) = 0.$

\item[{\rm (d)}] ${\mathcal S} = \{ 0\}$ if and only if
$\max_{X\in {\mathcal S}} r(X) = 0.$

\item[{\rm (e)}] ${\mathcal H}$ has a matrix $X \succ 0$  $(X \prec 0)$ if and only if
$\max_{X\in {\mathcal H}} i_{+}(X) = m  \ \left(\max_{X\in
{\mathcal H}} i_{-}(X) = m \right)\!.$

\item[{\rm (f)}] All $X\in {\mathcal H}$ satisfy $X \succ0$ $(X \prec 0)$ if and only if
$\min_{X\in {\mathcal H}} i_{+}(X) = m \ \left(\min_{X\in {\mathcal H}}
i_{-}(X) = m\, \right)\!.$

\item[{\rm (g)}] ${\mathcal H}$ has a matrix  $X \succcurlyeq 0$ $(X \preccurlyeq 0)$ if and only if
$\min_{X \in {\mathcal H}} i_{-}(X) = 0 \ \left(\min_{X\in
{\mathcal H}} i_{+}(X) = 0 \,\right)\!.$

\item[{\rm (h)}] All $X\in {\mathcal H}$ satisfy $X \succcurlyeq 0$ $(X \preccurlyeq 0)$  if and only if
$\max_{X\in {\mathcal H}} i_{-}(X) = 0 \ \left(\max_{X\in {\mathcal H}}
i_{+}(\,X) = 0\, \right)\!.$

\item[{\rm (i)}] The following hold
\begin{align}
& {\cal S}_1 \cap {\cal S}_2 \neq \emptyset  \ \ \Leftrightarrow  \ \
\min_{X_1\in {\cal S}_1, \, X_2\in {\cal S}_2} r(\,X_1 - X_2\,) = 0,
\label{115}
\\
& {\cal S}_1 \subseteq {\cal S}_2  \ \ \Leftrightarrow
\ \ \max_{X_1\in {\cal S}_1} \min_{X_2\in {\cal S}_2} r(\,X_1 - X_2\,) = 0,
\label{116}
\\
& {\cal S}_1 \supseteq {\cal S}_2 \Leftrightarrow \ \ \max_{X_2\in {\cal S}_2}
\min_{X_1\in {\cal S}_1} r(\,X_1 - X_2\,) = 0,
\label{117}
\\
& there \ exist \ X_1 \in {\cal S}_1 \ and \ X_2 \in {\cal S}_2 \ such \ that \ X_1 \succ X_2  \Leftrightarrow \!
\max_{X_1\in {\cal S}_1, \, X_2\in {\cal S}_2} \!\!\!i_{+}(\,X_1 - X_2\,) = m,
\label{118}
\\
& there \ exist \ X_1 \in {\cal S}_1 \ and \ X_2 \in {\cal S}_2 \ such \ that \ X_1 \succcurlyeq X_2 \Leftrightarrow\!
\min_{X_1\in {\cal S}_1, \, X_2\in {\cal S}_2}\!\!\!i_{-}(\,X_1 - X_2\,) = 0. \
\label{119}
\end{align}

\end{enumerate}
\end{lemma}

\begin{lemma} \label{T17}
Let $A\in {\mathbb C}^m_{{\rm H}},$  $B \in {\mathbb C}^n_{{\rm
H}},$ $Q \in {\mathbb C}^{m \times n},$ and assume that $P\in
{\mathbb C}^{m\times m}$ is nonsingular$.$ Then$,$
\begin{align}
& i_{\pm}(PAP^{*}) =  i_{\pm}(A),
\label{120}
\\
& i_{\pm}(\lambda A)  = \left\{ \begin{array}{ll}  i_{\pm}(A)  &  if  \ \lambda  >0
 \\ i_{\mp}(A)  & if  \ \lambda < 0
 \end{array},
\right.
\label{121}
\\
& i_{\pm}\!\left[ \begin{array}{cc}  A  & 0 \\  0  & B \end{array}
\right] = i_{\pm}(A) + i_{\pm}(B),
\label{122}
\\
& i_{+}\!\left[ \begin{array}{cc}  0  & Q \\  Q^{*} & 0 \end{array}
\right] =  i_{-}\!\left[ \begin{array}{cc}  0  & Q \\  Q^{*} & 0 \end{array}
\right] = r(Q).
\label{123}
\end{align}
\end{lemma}

\begin{lemma} [\cite{MS}] \label{T18}
Let $A \in {\mathbb C}^{m \times n},$ $B\in {\mathbb C}^{m \times k},$
$C \in {\mathbb C}^{l \times n}$ and  $D \in {\mathbb C}^{l \times k}.$ Then$,$
\begin{align}
 r[\, A, \, B \,] & = r(A)+ r(E_AB) = r(B) + r(E_BA),
\label{124}
\\
r\!\left[\!\! \begin{array}{cc} A \\ C \end{array} \!\!\right]
& = r(A) + r(CF_A)  =  r(C) + r(AF_C),
\label{125}
\\
r\!\left[\!\! \begin{array}{cc} A & B \\ C & 0 \end{array} \!\!\right]
& = r(B) + r(C) + r(E_BAF_C).
\label{126}
\end{align}
\end{lemma}

\begin{lemma} [\cite{T-laa10}] \label{T19}
Let $A \in {\mathbb C}_{{\rm H}}^{m},$ $B \in {\mathbb C}^{m\times n},$
 $D \in {\mathbb C}_{{\rm H}}^{n},$  and let
\begin{align}
M_1 = \left[\!\!\begin{array}{cc}  A  & B  \\ B^{*}  & 0 \end{array}
\!\!\right]\!, \ \ M_2 = \left[\!\!\begin{array}{cc}  A  & B  \\ B^{*}  & D \end{array}
\!\!\right]\!.
\label{127}
\end{align}
Then$,$ the following expansion formulas hold
\begin{align}
 i_{\pm}(M_1) = r(B) + i_{\pm}(E_BAE_B),   \ \ \ \ \ \ \ \ \ \ \ \ \ \ r(M_1) = 2r(B) + r(E_BAE_B), \ \ \ \ \ \  \ \ \ \
\label{128}
\\
i_{\pm}(M_2) = i_{\pm}(A) + i_{\pm}\!\left[\!\!\begin{array}{cc} 0 & \!\! E_AB
 \\
 B^{*}E_A & \!\!  D - B^{*}A^{\dag}B \end{array}\!\!\right]\!, \ \
 r(M_2)  = r(A) + r\!\left[\!\!\begin{array}{cc} 0 & \!\! E_AB
 \\
 B^{*}E_A & \!\! D - B^{*}A^{\dag}B  \end{array}\!\!\right]\!.
\label{129}
\end{align}
Under the condition $A \succcurlyeq 0,$
\begin{align}
i_{+}(M_1) = r[\, A, \, B \,],  \ \ i_{-}(M_1) = r(B), \ \ r(M_1) = r[\, A, \, B \,] +  r(B).
\label{130}
\end{align}
Under the condition $\R(B) \subseteq \R(A),$
\begin{align}
i_{\pm}(M_2) = i_{\pm}(A) + i_{\pm}(\, D - B^{*}A^{\dag}B \,), \ \ r(M_2)  = r(A) + r(\, D - B^{*}A^{\dag}B \,).
\label{131}
\end{align}
\end{lemma}

Some general rank and inertia expansion formulas derived from
(\ref{124})--(\ref{129}) are given below
\begin{align}
& r\!\left[\!\begin{array}{cc}
A & B\\ E_PC & 0 \end{array} \!\right]  = r\!\left[\!\!\begin{array}{ccc} A & B & 0  \\ C & 0 & P
\end{array} \!\!\right] - r(P),
\label{132}
\\
& r\!\left[\!\!\begin{array}{cc} A & BF_Q \\ C & 0 \end{array} \!\!\right] = r\!\left[\!\!
\begin{array}{ccc} A & B \\ C & 0 \\ 0 & Q
\end{array} \!\!\right] - r(Q),
\label{133}
\\
& r\!\left[\!\!\begin{array}{cc}
A & BF_Q \\ E_PC & 0 \end{array} \!\!\right] = r\!\left[\!\! \begin{array}{ccc} A & B & 0  \\ C & 0 & P \\ 0 & Q & 0
\end{array} \!\!\right] - r(P) - r(Q),
\label{134}
\\
& i_{\pm}\!\left[\!\! \begin{array}{cc} A
 & BF_P \\  F_P B^*  & 0 \end{array} \!\!\right]  = i_{\pm}\!\left[\!\! \begin{array}{ccc} A
 & B & 0 \\  B^* & 0 & P^* \\  0 & P & 0 \end{array} \!\!\right] - r(P),
\label{135}
\\
& i_{\pm}\!\left[\!\! \begin{array}{cc} E_{Q}AE_{Q}
 & E_{Q}B \\  B^*E_Q & D \end{array} \!\!\right]  = i_{\pm}\!\left[\!\! \begin{array}{ccc} A
 & B & Q \\  B^* & D & 0 \\  Q^* & 0 & 0 \end{array} \!\!\right] - r(Q).
\label{136}
\end{align}
We shall use them to simplify ranks and inertias of block  matrices involving Moore--Penrose inverses of matrices.

\begin{lemma} [\cite{KM}]\label{T110}
Let $A, \, B\in \mathbb C^{m\times n}$ be given$.$ Then$,$ the following hold$.$
\begin{enumerate}
\item[{\rm(a)}]  Eq. {\rm (\ref{11})} has a solution $X \in \mathbb C_{{\rm H}}^{n}$  if and only if ${\mathscr R}(B) \subseteq
{\mathscr R}(A)$ and $AB^{*} =BA^{*}.$ In this case$,$ the general
Hermitian solution of {\rm (\ref{11})} can be written as
\begin{equation}
X = A^{\dag}B + (A^{\dag}B)^{*} - A^{\dag}BA^{\dag}A + F_AUF_A,
\label{137}
\end{equation}
where $U\in \mathbb C_{{\rm H}}^{n}$ is arbitrary$.$

\item[{\rm(b)}] The matrix equation in {\rm (\ref{11})}
has a solution $0 \preccurlyeq X \in {\mathbb C}_{{\rm H}}^{n}$ if and only if ${\mathscr R}(B) \subseteq {\mathscr
R}(A),$ $AB^{*} \succcurlyeq 0$ and $\R(AB^{*}) = \R(BA^{*}) =\R(B).$ In this case$,$ the
general solution $0 \preccurlyeq X \in {\mathbb C}_{{\rm H}}^{n}$ of {\rm (\ref{11})} can be written as
\begin{equation}
X = B^{*}(AB^{*})^{\dag}B   + F_AUF_A,
\label{138}
\end{equation}
where $0 \preccurlyeq U\in {\mathbb C}_{{\rm H}}^{n}$ is arbitrary$.$

\end{enumerate}
\end{lemma}

\begin{lemma} \label{T111}
Eq. {\rm (\ref{12})}
has a solution $X \in {\mathbb C}_{{\rm H}}^{n}$  if and only if
 ${\mathscr R}(B) \subseteq {\mathscr R}(A),$ or equivalently$,$ $AA^{\dag}B = B.$
In this case$,$ the general Hermitian solution of $AXA^{*} = B$ can be written as
\begin{align}
X = A^{\dag}B(A^{\dag})^{*} +  F_AU + U^{*}F_A,
\label{139}
\end{align}
where  $U \in \mathbb C^{n \times n}$ is arbitrary$.$
\end{lemma}

\begin{lemma} [\cite{T-laa11}]\label{T112}
Let $A_j \in {\mathbb C}^{m_j\times n}$ and $B_j\in \mathbb C_{{\rm H}}^{m_j}$ be given$,$ $j =1, \,2,$  and
assume that
\begin{equation}
A_1X_1A_1^{*} = B_1 \ \    and  \ \ A_2X_2A_2^{*} = B_2
\label{140}
\end{equation}
are solvable for $X_1, \, X_2\in {\mathbb C}_{{\rm H}}^{n}.$ Also define
\begin{equation}
{\mathcal S}_j = \left\{\,  X_j\in {\mathbb C}_{{\rm H}}^{n}  \ | \ A_jX_jA_j^{*} = B_j \, \right\}, \ \  j =1, \, 2, \ \
M = \left[\!\!\begin{array}{ccccc}
B_1 & 0 &  A_1
\\
0 & -B_2  & A_2
\\
A_1^{*} & A^{*}_2 & 0\end{array}\!\!\right]\!.
\label{141}
\end{equation}
Then$,$
\begin{align}
\max_{X_1 \in {\mathcal S}_1,  X_2 \in {\mathcal S}_2}
r(\, X_1 - X_2\, ) & = \min\left\{ \, n, \ \  r(M) + 2n - 2r(A_1) - 2r( A_2) \, \right\},
\label{142}
\\
\min_{X_1 \in {\mathcal S}_1,  X_2 \in {\mathcal S}_2}
r(\,X_1 -X_2\,) & = r(M) - 2r[\, A^{*}_1, \,  A^{*}_2\,],
\label{143}
\\
\max_{X_1 \in {\mathcal S}_1,  X_2 \in {\mathcal S}_2}
i_{\pm}(\,X_1 -X_2\,) & = i_{\pm}(M) + n - r(A_1) - r( A_2),
\label{144}
\\
\min_{X_1 \in {\mathcal S}_1,  X_2 \in {\mathcal S}_2}
i_{\pm}(\,X_1 -X_2\,) & = i_{\pm}(M) - r[\, A^{*}_1, \,  A^{*}_2\,].
\label{145}
\end{align}
Consequently$,$ the following hold$.$
\begin{enumerate}
\item[{\rm(a)}] There exist  $X_1 \in {\mathcal S}_1$ and $X_2 \in {\mathcal S}_2$
 such that $X_1 -X_2$ is nonsingular if and only if $r(M) \geqslant 2r(A_1) + 2r( A_2) -n.$

\item[{\rm(b)}] $X_1 - X_2$ is nonsingular for all $X_1 \in {\mathcal S}_1$ and $X_2 \in {\mathcal S}_2$
if and only if $r(M) = 2r[\, A^{*}_1, \,  A^{*}_2\,] +n.$

\item[{\rm(c)}]  There exist  $X_1 \in {\mathcal S}_1$ and $X_2 \in {\mathcal S}_2$ such that
$X_1 = X_2$ if and only if
${\mathscr R}(B_j) \subseteq {\mathscr R}(A_j)$ and $r(M) = 2r[\, A^{*}_1, \,  A^{*}_2\,],$ $j =1, \,2.$

\item[{\rm(d)}] The rank of $X_1 - X_2$ is invariant for all $X_1 \in {\mathcal S}_1$ and $X_2 \in {\mathcal S}_2$
if and only if
$r(M) = 2r[\, A^{*}_1, \,  A^{*}_2\,] - n$ or $r(A_1) = r(A_2) = n.$

\item[{\rm(e)}] There exist  $X_1 \in {\mathcal S}_1$ and $X_2 \in {\mathcal S}_2$
such that $X_1 \succ X_2$ $(X_1 \prec X_2)$ if and only if $i_{+}(M) = r(A_1) + r( A_2)$\,
$\left(\, i_{-}(M) = r(A_1) + r( A_2)\,\right).$

\item[{\rm(f)}]  $X_1 \succ X_2$ $(X_1 \prec X_2)$ for all  $X_1 \in {\mathcal S}_1$ and $X_2 \in {\mathcal S}_2$
if and only if $i_{+}(M) = r[\, A^{*}_1, \,  A^{*}_2\,] +n$
$\left(\, i_{-}(M) = r[\, A^{*}_1, \,  A^{*}_2\,] + n \, \right).$

\item[{\rm(g)}] There exist  $X_1 \in {\mathcal S}_1$ and $X_2 \in {\mathcal S}_2$
 such that $X_1 \succcurlyeq X_2$ $(X_1 \preccurlyeq X_2)$ if and only if
 $i_{-}(M) = r[\, A^{*}_1, \,  A^{*}_2\,]$ $\left(\,i_{+}(M) =
 r[\, A^{*}_1, \, A^{*}_2\,]\,\right).$

\item[{\rm(h)}] $X_1 \succcurlyeq X_2$ $(X_1 \preccurlyeq X_2)$ for all  $X_1 \in {\mathcal S}_1$ and $X_2 \in {\mathcal S}_2$   if and only if $i_{-}(M) = r(A_1) + r( A_2) -n$ $\left( \, i_{+}(M) =  r(A_1) + r( A_2) -n \, \right).$

\item[{\rm (i)}] $i_{+}(\,X_1 - X_2\,)$ is invariant for all  $X_1 \in {\mathcal S}_1$ and $X_2 \in {\mathcal S}_2$    $\Leftrightarrow$ $i_{-}(\,X_1 - X_2\,)$ is invariant for all  $X_1 \in {\mathcal S}_1$ and $X_2 \in {\mathcal S}_2$  $\Leftrightarrow$ $r(A_1)  = r(A_2) = n.$
\end{enumerate}
\end{lemma}

\begin{lemma}  \label{T113}
 Let $A \in {\mathbb C}^m_{{\rm H}}$ and $B\in \mathbb C^{m\times n}$ be given$,$ and denote
 $M = \left[\!\! \begin{array}{cc}  A & B \\ B^* & 0\end{array} \!\!\right]\!.$ Then$,$ the following hold$.$
\begin{enumerate}
\item[{\rm(a)}]  {\rm \cite{T-laa10,TL}} The extremal ranks and inertias of $A - BXB^*$ subject to
$X \in {\mathbb C}^n_{{\rm H}}$  are given by
  \begin{align}
\max_{X\in {\mathbb C}^n_{{\rm H}}} r(\,A-BXB^*\,) & = r[\, A, \,B \,],
 \label{146}
\\
\min_{X \in {\mathbb C}^n_{{\rm H}}} r(\,A-BXB^*\,) & = 2r[\, A, \, B\,]- r(M),
\label{147}
 \\
\max_{X\in {\mathbb C}^n_{{\rm H}}}i_{\pm}(\,A-BXB^*\,)& = i_{\pm}(M),
\label{148}
 \\
\min_{X \in {\mathbb C}^n_{{\rm H}}}  i_{\pm}(\,A-BXB^*\,) & =  r[\, A, \, B \,]  -
 i_{\mp}(M).
 \label{149}
\end{align}

\item[{\rm(b)}]  {\rm \cite{T-mcm}} The extremal ranks and inertias of $A \pm BXB^*$ subject to
$0 \preccurlyeq X \in {\mathbb C}^n_{{\rm H}}$  are given by
\begin{align}
&\!\!\!\!\!\!\max_{0 \preccurlyeq X\in {\mathbb C}^{n}_{{\rm H}}}\!\!\!\!\!r(\, A + BXB^* \,)
= r[\, A, \, B\,], \  \min_{0 \preccurlyeq X\in {\mathbb C}^{n}_{{\rm H}}}\!\!\!\!r(\, A + BXB^* \,)
 = i_{+}(A) + r[\, A, \, B\,] - i_{+}(M),
\label{150}
\\
&\!\!\!\!\!\!\max_{0 \preccurlyeq X\in {\mathbb C}^{n}_{{\rm H}}}\!\!\!\!i_{+}(\, A + BXB^* \,)  = i_{+}(M), \
  \min_{0 \preccurlyeq X\in {\mathbb C}^{n}_{{\rm H}}}\!\!\!i_{+}(\, A + BXB^* \,)  = i_{+}(A),
\label{151}
\\
&\!\!\!\!\!\!\max_{0 \preccurlyeq X\in {\mathbb C}^{n}_{{\rm H}}}\!\!\!\!i_{-}(\, A + BXB^* \,)  = i_{-}(A), \
\min_{0 \preccurlyeq X\in {\mathbb C}^{n}_{{\rm H}}}\!\!\!i_{-}(\, A + BXB^* \,)  = r[\, A, \, B\,] - i_{+}(M),
 \label{152}
\\
&\!\!\!\!\!\!\max_{0 \preccurlyeq X\in {\mathbb C}^{n}_{{\rm H}}}\!\!\!\!r(\, A - BXB^* \,)
= r[\, A, \, B\,],  \  \min_{0 \preccurlyeq X\in {\mathbb C}^{n}_{{\rm H}}}\!\!\!\!r(\, A - BXB^* \,) =
i_{-}(A) + r[\, A, \, B\,] - i_{-}(M),
\label{153}
\\
&\!\!\!\!\!\!\max_{0 \preccurlyeq X\in {\mathbb C}^{n}_{{\rm H}}}\!\!\!\!i_{+}(\, A - BXB^* \,)  = i_{+}(A), \
\min_{0 \preccurlyeq X\in {\mathbb C}^{n}_{{\rm H}}}\!\!\!i_{+}(\, A - BXB^* \,) =  r[\,A, \, B\,] -  i_{-}(M),
\label{154}
\\
&\!\!\!\!\!\!\max_{0 \preccurlyeq X\in {\mathbb C}^{n}_{{\rm H}}}\!\!\!\!i_{-}(\, A - BXB^* \,)  = i_{-}(M),  \
\min_{0 \preccurlyeq X\in {\mathbb C}^{n}_{{\rm H}}}\!\!\!i_{-}(\, A -
 BXB^* \,) = i_{-}(A).
 \label{155}
\end{align}
\end{enumerate}
\end{lemma}

\begin{lemma}  [\cite{LT-nla11,T-laa10}] \label{T114}
Let $A \in {\mathbb C}_{{\rm H}}^{m}$, $B \in \mathbb C^{m\times n}$ and $C \in \mathbb C^{m\times k}$
 be given$.$ Then$,$
\begin{align}
  \max_{X \in {\mathbb C}_{{\rm H}}^{n}, \
Y \in {\mathbb C}_{{\rm H}}^{k}} \!\!\!r(\, A - BXB^* - CYC^* \,) & = r[\,A,\, B,\, C \,],
\label{156}
\\
 \min_{X \in {\mathbb C}_{{\rm H}}^{n}, \ Y \in {\mathbb C}_{{\rm H}}^{k}}\!\!\!r(\, A - BXB^* - CYC^* \,)
&  = 2r[\,A,\, B,\, C \,] +  r\!\left[\!\!\begin{array}{cc} A & B \\ C^* & 0
  \end{array}\!\!\right] - r\!\left[\!\!\begin{array}{ccc} A & B & C \\ B^* &
  0 & 0
  \end{array}\!\!\right] \nb
  \\
& \ \ \ \  - r\!\left[\!\!\begin{array}{ccc} A & B & C \\ C^* &
  0 & 0
  \end{array}\!\!\right]\!,
 \label{157}
  \\
 \max_{X \in {\mathbb C}_{{\rm H}}^{n}, \ Y \in {\mathbb C}_{{\rm H}}^{k}} \!\!\!i_{\pm}(\, A - BXB^* - CYC^* \,)
& =  i_{\pm}\!\left[\!\!\begin{array}{ccc} A & B & C \\ B^* & 0 & 0 \\ C^* & 0 & 0
  \end{array}\!\!\right]\!,
\label{158}
\\
 \min_{X \in {\mathbb C}_{{\rm H}}^{n}, \ Y \in {\mathbb C}_{{\rm H}}^{k}} \!\!\!i_{\pm}(\, A - BXB^* - CYC^* \,)
 & =  r[\,A,\, B,\, C \,] - i_{\mp}\!\left[\!\!\begin{array}{ccc} A & B & C \\ B^* & 0 & 0 \\ C^* & 0 & 0
  \end{array}\!\right]\!.
\label{159}
\end{align}
\end{lemma}

\begin{lemma}[\cite{LT-jota,T-laa11}]\label{T115}
Let $A \in {\mathbb C}_{{\rm H}}^{m}$, $B \in {\mathbb C}^{m\times n}$ and $C \in {\mathbb C}^{p \times m}$
be given and assume that ${\mathscr R}(B)\subseteq {\mathscr R}(C^*).$ Then$,$
\begin{align}
\max_{X\in {\mathbb C}^{n \times p}}r[\,A-BXC-(BXC)^{*}\,] & = \min \left\{ r[\,A,\, C^{*}\,],
 \ \ r\!\left[\!\!\begin{array}{cc} A & B
\\B^{*} & 0
\end{array}\!\!\right] \right\}\!,
\label{160}
\\
\min_{X\in {\mathbb C}^{n \times p}}r[\,A-BXC-(BXC)^{*}\,] & =
2r[\,A,\, C^{*}\,] + r\!\left[\!\!\begin{array}{cc} A & B
\\B^{*} & 0
\end{array}\!\!\right] - 2r\!\left[\!\!\begin{array}{cc} A & B
\\C & 0
\end{array}\!\!\right]\!,
\label{161}
\\
\max_{X\in {\mathbb C}^{n \times p}}i_{\pm}[\,A-BXC-(BXC)^{*}\,] &
= i_{\pm}\!\left[\!\!\begin{array}{ccc} A
 & B  \\  B^{*}  & 0
   \end{array}\!\!\right]\!,
\label{162}
\\
\min_{X\in {\mathbb C}^{n \times p}}i_{\pm}[\,A-BXC-(BXC)^{*}\,]&
= r[\,A,\, C^{*}\,] + i_{\pm}\!\left[\!\!\begin{array}{ccc} A
 & B  \\  B^{*}  & 0
   \end{array}\!\!\right]- r\!\left[\!\!\begin{array}{cc} A & B
\\ C & 0
\end{array}\!\!\right]\!,
\label{163}
\end{align}
and
\begin{align}
\max_{X\in {\mathbb C}^{n \times m}}r[\,A - BX - (BX)^{*}\,] & = \min \left\{m,
 \ \ r\!\left[\!\!\begin{array}{cc} A & B
\\B^{*} & 0
\end{array}\!\!\right] \right\}\!,
\label{164}
\\
\min_{X\in {\mathbb C}^{n \times m}}r[\,A - BX - (BX)^{*}\,] & = r\!\left[\!\!\begin{array}{cc} A & B
\\B^{*} & 0
\end{array}\!\!\right] - 2r(B),
\label{165}
\\
\max_{X\in {\mathbb C}^{n \times m}}i_{\pm}[\,A  - BX - (BX)^{*}\,] &
= i_{\pm}\!\left[\!\!\begin{array}{ccc} A & B  \\  B^{*}  & 0
   \end{array}\!\!\right]\!,
\label{166}
\\
\min_{X\in {\mathbb C}^{n \times m}}i_{\pm}[\,A -  BX - (BX)^{*}\,]&
= i_{\pm}\!\left[\!\!\begin{array}{ccc} A
 & B  \\  B^{*}  & 0
   \end{array}\!\!\right]- r(B).
\label{167}
\end{align}
\end{lemma}

\section{Properties of Hermitian solutions and  Hermitian definite solutions of $AX= B$}
\renewcommand{\theequation}{\thesection.\arabic{equation}}
\setcounter{section}{2} \setcounter{equation}{0}

 Some formulas for calculating the ranks and inertias of  Hermitian solutions and Hermitian definite solutions of the matrix equation in (\ref{11}) were established in \cite{T-laa11}. In this section, we reconsider the ranks and inertias of these
 solutions  and give a group of complete results.

\begin{theorem} \label{T21}
Assume that {\rm (\ref{11})} has a Hermitian  solution$,$
and let $P\in {\mathbb C}^{n}_{{\rm H}}.$ Also$,$ define
\begin{align}
{\mathcal S} = \{\, X\in {\mathbb C}^{n}_{{\rm H}} \ | \  AX = B \, \}.
\label{21}
\end{align}
Then$,$
\begin{align}
\max_{X \in {\mathcal S}}r(\, X - P\,) & = r(\,B -  AP \,) - r(A) + n,
\label{22}
\\
\min_{X \in {\mathcal S}}r(\, X - P \,) & =  2r(\, B - AP \,) - r(\,BA^* - APA^* \,),
\label{23}
\\
\max_{X \in {\mathcal S}} i_{\pm}(\, X - P \,) & = i_{\pm}(\,BA^* - APA^* \,) - r(A) + n,
\label{24}
\\
\min_{X \in {\mathcal S}} i_{\pm}(\, X - P \,) & =  r(\,B -  AP \,) - i_{\mp}(\,BA^* - APA^* \,).
\label{25}
\end{align}
In consequence$,$ the following hold$.$
\begin{enumerate}
\item[{\rm(a)}] There exists an $X \in {\mathcal S}$ such that $X - P$ is nonsingular
   if and only if $\R(\, AP - B\,) = \R(A).$

\item[{\rm(b)}]  $X - P$ is nonsingular for all $X \in {\mathcal S}$ if and only if
$2r(\,B -  AP\,) = r(\,BA^* - APA^* \,) + n.$

\item[{\rm(c)}] There exists an  $X \in {\mathcal S}$ such that $X \succ P$ $(X \prec P)$
   if and only if
$$
\R(\, BA^* - APA^* \,) = \R(A) \ \ and \ \ BA^* \succcurlyeq APA^*
$$
$$
\left( \, \R(\, BA^* - APA^* \,) = \R(A) \ \ and \ \ BA^* \preccurlyeq APA^* \,\right).
$$

\item[{\rm(d)}]  $X \succ P$ $(X \prec P)$ holds for all $X \in {\mathcal S}$ if and only if
$$
r(\,B -  AP\,) = n  \ \ and \ \  BA^* \succcurlyeq APA^*\ \  \left(\,r(\, B - AP\,) = n  \ \ and \ \
 AB^* \preccurlyeq APA^* \, \right).
$$

\item[{\rm(e)}] There exists an $X \in {\mathcal S}$ such that $X \succcurlyeq P$
 $(X \preccurlyeq P)$ if and only if
$$
 \R(\,B -  AP \,) = \R(\, BA^* - APA^* \,)  \ \ and  \ \
  BA^* \succcurlyeq APA^*
$$
$$
 \left( \, \R(\, B - AP \,) = \R(\,BA^* - APA^*\,)  \ \ and  \ \
  BA^* \preccurlyeq APA^*  \, \right)\!.
$$

\item[{\rm(f)}]  $X \succcurlyeq P$ $(X \preccurlyeq P)$  holds for all $X \in {\mathcal S}$ if and only if
$$
 BA^* \succcurlyeq APA^* \ \ and  \ \   r(A) = n  \ \ \left(\, BA^* \preccurlyeq APA^* \ \ and  \ \   r(A) = n  \, \right).
$$
\end{enumerate}
In particular$,$ the following hold$.$
\begin{enumerate}
\item[{\rm(g)}] There exists an $X \in {\mathcal S}$ such that $X$ is nonsingular if and only if
$\R(B) = \R(A).$

\item[{\rm(h)}]  $X$ is nonsingular for all $X \in {\mathcal S}$ if and only if $r(\,B\,)= n.$

\item[{\rm(i)}] There exists an $X \in {\mathcal S}$ such that $X \succ 0$ $(X \prec 0)$  if and only if
$$
\R(B) = \R(A) \ \ and \ \ BA^* \succcurlyeq 0 \ \
\left( \, \R(B) = \R(A) \ \ and \ \ BA^* \preccurlyeq 0 \,\right).
$$

\item[{\rm(j)}]  $X \succ 0$ $(X \prec 0)$ holds for all $X \in {\mathcal S}$ if and only if
$$
r(B) = n  \ \ and \ \  BA^* \succcurlyeq 0 \ \  \left(\,r(B) = n  \ \ and \ \
 AB^* \preccurlyeq 0\, \right).
$$

\item[{\rm(k)}] There exists an $X \in {\mathcal S}$ such that $X \succcurlyeq 0$ $(X \preccurlyeq 0)$
 if and only if
$$
 \R(B) = \R(BA^*)  \ \ and  \ \  BA^* \succcurlyeq 0\ \  \left( \, \R(B) = \R(BA^*)  \ \ and  \ \
  BA^* \preccurlyeq 0 \, \right)\!.
$$

\item[{\rm(l)}]  $X \succcurlyeq 0$ $(X \preccurlyeq 0)$  holds for all $X \in {\mathcal S}$ if and only if
$$
 BA^* \succcurlyeq 0 \ \ and  \ \   r(A) = n  \ \ \left(\, BA^* \preccurlyeq 0 \ \ and  \ \   r(A) = n  \, \right).
$$
\end{enumerate}
\end{theorem}

\noindent {\bf Proof} \ By Lemma \ref{T110}(a), $X - P$ can be written as
\begin{align}
X - P   =  X_0 - P + F_AUF_A,
\label{26}
\end{align}
where  $X_0=  A^{\dag}B + (A^{\dag}B)^{*} - A^{\dag}BA^{\dag}A$ and
$U \in {\mathbb C}_{{\rm H}}^{n}$ is  arbitrary.
Applying Lemma \ref{T113}(a) to (\ref{26}) gives
\begin{align}
\max_{X \in {\mathcal S}}r(\, X - P \,)
& = \max_{U \in {\mathbb C}_{{\rm H}}^{n}}r\!\left(\, X_0 - P + F_AUF_A \right) = r[\, X_0 - P,\, F_A \,],
\label{27}
\\
\min_{X \in {\mathcal S}}r(\, X - P \,)
& = \min_{U \in {\mathbb C}_{{\rm H}}^{n}}r\!\left(\, X_0 - P + F_AUF_A \right) =
 2r[\, X_0 - P,\, F_A \,] -  r\!\left[\!\! \begin{array}{ccc}
X_0 - P  & F_A
\\
F_A  & 0
\end{array}\!\!\right]\!,
\label{28}
\\
\max_{X \in {\mathcal S}}i_{\pm}(\, X - P \,)
 & = \min_{U \in {\mathbb C}_{{\rm H}}^{n}}i_{\pm}\!\left(\, X_0 - P + F_AUF_A \right)
 = i_{\pm}\!\left[\!\! \begin{array}{cc}
X_0 - P  & F_A
\\
F_A  & 0
\end{array}\!\!\right]\!,
\label{29}
\\
\min_{X \in {\mathcal S}}i_{\pm}(\, X - P \,)
& = \min_{U \in {\mathbb C}_{{\rm H}}^{n}}i_{\pm}\!\left(\, X_0 - P + F_AUF_A \right) =
 r[\, X_0 - P,\, F_A \,] -  i_{\mp}\!\left[\!\! \begin{array}{ccc}
X_0 - P  & F_A
\\
F_A  & 0
\end{array}\!\!\right]\!.
\label{210}
\end{align}
Applying (\ref{124}) and (\ref{128}) to the block matrices in (\ref{27})--(\ref{210}) and simplifying by
(\ref{120}) and elementary block matrix operations, we obtain
\begin{align}
 r[\, X_0 - P,\, F_A \,] & = r(\, A^{\dag}AX_0 - A^{\dag}AP\,) + r(F_A)  =
 r(\, B - AP\,) + n - r(A),
\label{211}
 \\
i_{\pm}\!\left[\!\! \begin{array}{ccc}
X_0 - P  & F_A
\\
F_A  & 0
\end{array}\!\!\right] & =  r(F_A) + i_{\pm}[\, A^{\dag}A(\, X_0 - P \,)A^{\dag}A\,] =
n -  r(A) + i_{\pm}(\,BA^* - APA^* \,),
\label{212}
\\
r\!\left[\!\! \begin{array}{ccc}
X_0 - P  & F_A
\\
F_A  & 0
\end{array}\!\!\right] & =  2r(F_A) + r[\, A^{\dag}A(\, X_0 - P \,)A^{\dag}A\,] =
2n -  2r(A) + r(\,BA^* - APA^* \,).
\label{213}
\end{align}
Substituting (\ref{211})--(\ref{213}) into (\ref{27})--(\ref{210}) and simplifying leads to (\ref{22})--(\ref{25}).
Results (a)--(l) follow from applying Lemma \ref{T16} to   (\ref{22})--(\ref{25}).
\qquad $\Box$

\begin{theorem} \label{T22}
Assume that {\rm (\ref{11})} has a Hermitian solution $X  \succcurlyeq 0,$
and let $ 0 \preccurlyeq  P\in {\mathbb C}^{n}_{{\rm H}}.$ Also$,$ define
\begin{align}
{\mathcal S} = \{\, 0 \preccurlyeq X\in {\mathbb C}^{n}_{{\rm H}} \ | \  AX = B \, \}, \ \ M = \left[\! \begin{array}{ccc}
 AB^*  & B   \\  B^*  & P  \end{array} \!\right]\!.
\label{214}
\end{align}
Then$,$
\begin{align}
\max_{X \in {\mathcal S}}r(\, X - P) & = r(\,B -  AP \,) - r(A) + n,
\label{215}
\\
\min_{X \in {\mathcal S}}r(\, X - P \,) & = i_{-}(M) +  r(\,B -  AP \,) - i_{+}(\,BA^* - APA^* \,),
\label{216}
\\
\max_{X \in {\mathcal S}} i_{+}(\, X - P \,) &  = i_{+}(\,BA^* - APA^* \,) - r(A) + n,
\label{217}
\\
\min_{X \in {\mathcal S}} i_{+}(\, X - P \,) & = i_{-}(M),
\label{218}
\\
\max_{X \in {\mathcal S}} i_{-}(\, X - P \,) &  = i_{+}(M) - r(B),
\label{219}
\\
\min_{X \in {\mathcal S}} i_{-}(\, X - P \,) & =  r(\,B -  AP \,) - i_{+}(\,BA^* - APA^* \,).
\label{220}
\end{align}
Consequently$,$ the following hold$.$
\begin{enumerate}
\item[{\rm(a)}] There exists an $X \in {\mathcal S}$ such that $X - P$ is nonsingular
   if and only if $\R(\, B - AP \,) = \R(A).$

\item[{\rm(b)}]  $X - P$ is nonsingular for all $X \in {\mathcal S}$ if and only if
$i_{-}(M) +  r(\,B -  AP \,) = i_{+}(\,BA^* - APA^* \,) + n.$

\item[{\rm(c)}] There exists an $X \in {\mathcal S}$ such that $X \succ P$
   if and only if $\R(\, BA^* - APA^* \,) = \R(A)$  and $BA^* \succcurlyeq APA^*.$

\item[{\rm(d)}]  $X \succ P$ holds for all $X \in {\mathcal S}$ if and only if
$i_{-}(M) = n.$

\item[{\rm(e)}] There exists an $X \in {\mathcal S}$ such that $X \prec P$
   if and only if $i_{+}(M) = r(B) + n.$

\item[{\rm(f)}]  $X \prec P$ holds for all $X \in {\mathcal S}$ if and only if
$r(\, B - AP \,) = n$  and $BA^* \preccurlyeq APA^*.$

\item[{\rm(g)}] There exists an $X \in {\mathcal S}$ such that $X \succcurlyeq P$  if and only if
$\R(\,B -  AP \,) = \R(\, BA^* - APA^* \,)$ and  $ BA^* \succcurlyeq APA^*.$

\item[{\rm(h)}]  $X \succcurlyeq P$  holds for all $X \in {\mathcal S}$ if and only if
$i_{+}(M) = r(B).$

\item[{\rm(i)}] There exists an $X \in {\mathcal S}$ such that $X \preccurlyeq P$  if and only if $M\succcurlyeq 0.$

\item[{\rm(j)}] $X \preccurlyeq P$  holds for all $X \in {\mathcal S}$ if and only if
$i_{+}(\,BA^* - APA^* \,) = n  - r(A).$
\end{enumerate}
In particular$,$ the following hold$.$
\begin{enumerate}
\item[{\rm(k)}] There exists an $X \in {\mathcal S}$ such that $X - I_n$ is nonsingular
   if and only if $\R(\, B - A \,) = \R(A).$

\item[{\rm(l)}]  $X - I_n$ is nonsingular for all $X \in {\mathcal S}$ if and only if
$i_{-}(\,BA^* - BB^* \,)  + r(\,B -  A \,) = i_{+}(\,BA^* - AA^* \,) + n.$

\item[{\rm(m)}] There exists an $X \in {\mathcal S}$ such that $X \succ I_n$
   if and only if $\R(\, BA^* - AA^* \,) = \R(A)$ and $BA^* \succcurlyeq AA^*.$

\item[{\rm(n)}]  $X \succ I_n$ holds for all $X \in {\mathcal S}$ if and only if
$BA^* \prec BB^*.$

\item[{\rm(o)}] There exists an $X \in {\mathcal S}$ such that $X \prec I_n$
   if and only if $\R(\,BA^* - BB^* \,)  = \R(B)$  and $BA^* \succcurlyeq BB^*.$

\item[{\rm(p)}]  $X \prec I_n$ holds for all $X \in {\mathcal S}$ if and only if
$r(\, B - A \,) = n$  and $BA^* \preccurlyeq AA^*.$

\item[{\rm(q)}] There exists an $X \in {\mathcal S}$ such that $X \succcurlyeq I_n$  if and only if
$\R(\,B -  A \,) = \R(\, BA^* - AA^* \,)$ and $ BA^* \succcurlyeq AA^*.$

\item[{\rm(r)}]  $X \succcurlyeq I_n$  holds for all $X \in {\mathcal S}$ if and only if
$r(B) = n $ and  $ BA^* \preccurlyeq BB^*.$

\item[{\rm(s)}] There exists an $X \in {\mathcal S}$ such that $X \preccurlyeq I_n$  if and only if
$BA^* \succcurlyeq BB^*.$

\item[{\rm(t)}] $X \preccurlyeq I_n$  holds for all $X \in {\mathcal S}$ if and only if $i_{+}(\,BA^* - AA^* \,) = n  - r(A).$
\end{enumerate}
\end{theorem}

\noindent {\bf Proof} \ By Lemma \ref{T110}(b), $X - P$ can be written as
\begin{align}
X - P   =  X_0 - P + F_AUF_A,
\label{221}
\end{align}
where  $X_0=  B^*(AB^*)^{\dag}B$ and $ 0 \preccurlyeq U \in {\mathbb C}_{{\rm H}}^{n}$ is  arbitrary.
Applying Lemma \ref{T113}(b) to (\ref{221}) gives
\begin{align}
\max_{X \in {\mathcal S}}r(\, X - P \,)
& = \max_{0 \preccurlyeq U \in {\mathbb C}_{{\rm H}}^{n}}r\!\left(\, X_0 - P + F_AUF_A \right) = r[\, X_0 - P,\, F_A \,],
\label{222}
\\
\min_{X \in {\mathcal S}}r(\, X - P \,)
& = \min_{0 \preccurlyeq U \in {\mathbb C}_{{\rm H}}^{n}}r\!\left(\, X_0 - P + F_AUF_A \right) =
 i_{+}(\, X_0 - P \,) +  r[\, X_0 - P,\, F_A \,] \nb
 \\
 & \ \ \  -  i_{+}\!\left[\!\! \begin{array}{ccc}
X_0 - P  & F_A
\\
F_A  & 0
\end{array}\!\!\right]\!,
\label{223}
\\
\max_{X \in {\mathcal S}}i_{+}(\, X - P \,)
 & = \min_{0 \preccurlyeq U \in {\mathbb C}_{{\rm H}}^{n}}i_{+}\!\left(\, X_0 - P + F_AUF_A \right)
 = i_{+}\!\left[\!\! \begin{array}{cc}
X_0 - P  & F_A
\\
F_A  & 0
\end{array}\!\!\right]\!,
\label{224}
\\
\min_{X \in {\mathcal S}}i_{+}(\, X - P \,)
 & = \min_{0 \preccurlyeq U \in {\mathbb C}_{{\rm H}}^{n}}i_{+}\!\left(\, X_0 - P + F_AUF_A \right)
 = i_{+}(\, X_0 - P \,),
\label{225}
\\
\max_{X \in {\mathcal S}}i_{-}(\, X - P \,)
 & = \min_{0 \preccurlyeq U \in {\mathbb C}_{{\rm H}}^{n}}i_{-}\!\left(\, X_0 - P + F_AUF_A \right)
 = i_{-}(\, X_0 - P \,),
\label{226}
\\
\min_{X \in {\mathcal S}}i_{-}(\, X - P \,)
& = \min_{0 \preccurlyeq U \in {\mathbb C}_{{\rm H}}^{n}}i_{\pm}\!\left(\, X_0 - P + F_AUF_A \right) =
 r[\, X_0 - P,\, F_A \,] -  i_{+}\!\left[\!\! \begin{array}{ccc}
X_0 - P  & F_A
\\
F_A  & 0
\end{array}\!\!\right]\!.
\label{227}
\end{align}
Applying (\ref{121}) and (\ref{131}) to $X_0 - P$ gives
\begin{align*}
i_{\pm}(\, X_0 - P \,)  =  i_{\mp}[\, P - B^*(AB^*)^{\dag}B^* \,]  =
i_{\mp}\!\left[\!\! \begin{array}{cc} AB^*  &  B^*   \\  B &  P \end{array} \!\!\right]-
i_{\mp}(AB^*) = i_{\mp}(M) - i_{\mp}(AB^*),
\end{align*}
so that
\begin{align}
i_{+}(\, X_0 - P \,) & = i_{-}(M), \ \ \ i_{+}(\, X_0 - P \,)  = i_{+}(M) - r(B).
\label{228}
\end{align}
Substituting (\ref{211})--(\ref{213}) and  (\ref{228}) into (\ref{222})--(\ref{227}) and simplifying leads
to (\ref{215})--(\ref{220}). Results (a)--(t) follow from applying Lemma \ref{T16} to   (\ref{215})--(\ref{220}).
\qquad $\Box$

\medskip

A general problem related to Hermitian solutions and Hermitian definite solutions of $AX= B$ is to establish
 formulas for calculating the extremal ranks and inertias of $P - QXQ^*$  subject to the Hermitian solutions and
 Hermitian definite solutions of $AX = B$. The results obtained can  be used to solve optimization problems of
 $P - QXQ^*$ subject to $AX= B$.

\section{Relations between Hermitian solutions of $AX= B$ and $CY = D$}
\renewcommand{\theequation}{\thesection.\arabic{equation}}
\setcounter{section}{3} \setcounter{equation}{0}

In order to compare Hermitian solutions of matrix equations, we first establish some
fundamental formulas  for calculating the extremal ranks and inertias of difference of Hermitian solutions
 of the two matrix equations $AX= B$ and $CY = D$, and then use them to characterize relationship
  between the Hermitian solutions.

\begin{theorem} \label{T31}
Assume that each of the matrix equations in {\rm (\ref{15})} has a Hermitian  solution$,$
and let
\begin{align}
{\mathcal S} = \{\, X\in {\mathbb C}^{n}_{{\rm H}} \ | \  AX = B \, \}, \ \
{\mathcal T} = \{ \, Y \in {\mathbb C}^{n}_{{\rm H}} \ | \  CY = D \, \}.
\label{31}
\end{align}
Also denote
$$
M = \left[\!\! \begin{array}{ccccc}
 AB^* & 0 & A
 \\
 0 & -CD^* & C
 \\
A^* &  C^* & 0
\end{array}\!\!\right]\!, \ \  N = \left[\!\! \begin{array}{cc}
A  & B
\\
C  & D
\end{array}\!\!\right]\!.
$$
Then$,$
\begin{align}
\max_{X \in {\mathcal S},\, Y \in {\mathcal T}}r(\, X - Y \,) & =  n + r(N) - r(A) - r(C),
\label{32}
\\
\min_{X \in {\mathcal S},\, Y \in {\mathcal T}}r(\, X - Y \,) & =
2r(N) + r(\, AD^* - BC^* \,) -  r\!\left[\!\! \begin{array}{cc}
A  & BA^*
\\
C  & DA^*
\end{array}\!\!\right] -  r\!\left[\!\! \begin{array}{ccc}
A  & BC^*
\\
C  & DC^*
\end{array}\!\!\right]\!,
\label{33}
\\
\max_{X \in {\mathcal S},\, Y \in {\mathcal T}} i_{\pm}(\, X - Y \,)
& =  n + i_{\pm}(M) - r(A) - r(C),
\label{34}
\\
\min_{X \in {\mathcal S},\, Y \in {\mathcal T}} i_{\pm}(\, X - Y \,)
& =  r(N) - i_{\mp}(M).
\label{35}
\end{align}
In consequence$,$ the following hold$.$
\begin{enumerate}
\item[{\rm(a)}] There exist $X \in {\mathcal S}$ and $Y\in {\mathcal T}$ such that $X - Y$ is nonsingular
   if and only if $r(N) =  r(A) + r(C).$

\item[{\rm(b)}]  $X - Y$ is nonsingular for all $X \in {\mathcal S}$ and $Y\in {\mathcal T}$ if and only if
$$
2r(N) + r(\, AD^* - BC^* \,)  =  r\!\left[\!\! \begin{array}{ccc}
A  & BA^*
\\
C  & DA^*
\end{array}\!\!\right] +  r\!\left[\!\! \begin{array}{ccc}
A  & BC^*
\\
C  & DC^*
\end{array}\!\!\right] + n.
$$

\item[{\rm(c)}] ${\mathcal S} \cap {\mathcal T} \neq \emptyset$ if and only if
$\R\!\left[\!\! \begin{array}{c}
B
\\
D
\end{array}\!\!\right] \subseteq \R\!\left[\!\! \begin{array}{ccc}
A
\\
C
\end{array}\!\!\right]\!, \ \  \left[\!\! \begin{array}{c}
A
\\
C
\end{array}\!\!\right]\![\, B^*, \ D^* \,] = \left[\!\! \begin{array}{c}
B
\\
D
\end{array}\!\!\right][\, A^*, \ C^* \,].$

\item[{\rm(d)}] There exist $X \in {\mathcal S}$ and $Y\in {\mathcal T}$ such that $X \succ Y$ $(X \prec Y)$
   if and only if $i_{+}(M) = r(A) + r(C)$  $\left( \, i_{-}(M) = r(A) + r(C)  \, \right).$

\item[{\rm(e)}]  $X \succ Y$ $(X \prec Y)$ holds for all $X \in {\mathcal S}$ and $Y\in {\mathcal T}$ if and only if
$i_{-}(M) = r(N) - n $ $\left(\, i_{+}(M) = r(N) - n \, \right)\!.$

\item[{\rm(f)}] There exist $X \in {\mathcal S}$ and $Y\in {\mathcal T}$ such that $X \succcurlyeq Y$
 $(X \preccurlyeq Y)$  holds if and only if
 $i_{+}(M) = r(N)$ $\left( \, i_{-}(M) = r(N) \, \right)\!.$

\item[{\rm(g)}]  $X \succcurlyeq Y$ $(X \preccurlyeq Y)$  holds for all
 $X \in {\mathcal S}$ and $Y\in {\mathcal T}$ if and only if
$i_{-}(M) = r(A)  + r(C) -n$ $\left(\, i_{+}(M)  = r(A)  + r(C) -n  \, \right).$
\end{enumerate}
\end{theorem}

\noindent {\bf Proof} \ By Lemma \ref{T110}(a), $X - Y$ can be written as
\begin{align}
X - Y   =  X_0 - Y_0 + F_AUF_A -  F_CVF_C,
\label{36}
\end{align}
where  $X_0= A^{\dag}B + (A^{\dag}B)^* - A^{\dag}BA^{\dag}A$ and
$Y_0 = C^{\dag}D + (C^{\dag}D)^* - C^{\dag}DC^{\dag}C,$
and $U, \ V \in {\mathbb C}_{{\rm H}}^{n}$ are arbitrary.
Applying Lemma \ref{T114} to (\ref{36}) gives
\begin{align}
 \max_{X \in {\mathcal S},\, Y \in {\mathcal T}}r(\, X - Y \,)
& = \max_{U, \, V}r\!\left(\, X_0 - Y_0 + F_AUF_A -  F_CVF_C \right) = r[\, X_0 - Y_0,\, F_A,\, F_C \,],
\label{37}
\\
 \min_{X \in {\mathcal S},\, Y \in {\mathcal T}}r(\, X - Y \,) & = \min_{U, \, V}r\!\left(\, X_0 - Y_0
 + F_AUF_A -  F_CVF_C \right) \nb
\\
& = 2r[\, X_0 - Y_0,\, F_A,\, F_C \,] + r\!\left[\!\! \begin{array}{cc}
 X_0 - Y_0   &  F_A
\\
F_C & 0
\end{array}\! \!\right]  -  r\!\left[\!\! \begin{array}{cccc}
 X_0 - Y_0   &  F_A & F_C
\\
F_A & 0 & 0
\end{array}\! \!\right] \nb
\\
& \ \ \ - r\!\left[\!\! \begin{array}{cccc}
 X_0 - Y_0   &  F_A & F_C
\\
F_C & 0 & 0
\end{array}\! \!\right]\!,
\label{38}
\\
 \max_{X \in {\mathcal S},\, Y \in {\mathcal T}}i_{\pm}(\, X - Y \,) & = \max_{U, \, V}i_{\pm}\!\left(\,X_0 - Y_0  + F_AUF_A -  F_CVF_C \right)  \nb
\\
&  = i_{\pm}\!\left[\!\! \begin{array}{cccc}
 X_0 - Y_0   &  F_A & F_C
\\
F_A & 0 & 0
\\
F_C & 0 & 0
\end{array}\! \!\right]\!,
\label{39}
\\
\min_{X \in {\mathcal S},\, Y \in {\mathcal T}}i_{\pm}(\, X - Y \,) &  =
\min_{U, \, V}i_{\pm}\!\left(\,X_0 - Y_0  + F_AUF_A -  F_CVF_C \right) \nb
\\
& = r[\,X_0 - Y_0 ,\, F_A,\, F_C \,] - i_{\mp}\!\left[\!\! \begin{array}{cccc}
 X_0 - Y_0   &  F_A & F_C
\\
F_A & 0 & 0
\\
F_C & 0 & 0
\end{array}\!\!\right]\!.
\label{310}
\end{align}
Applying (\ref{132})--(\ref{134}) to the block matrices in (\ref{37})--(\ref{310}) and simplifying
by Lemma \ref{T17}, (\ref{120}) and elementary block matrix operations, we obtain
\begin{align}
r[\,X_0 - Y_0,\, F_A,\, F_C \,] &= r\!\left[\!\! \begin{array}{ccc}
 X_0 - Y_0   &  I_n & I_n
\\
0 & A & 0
\\
0 & 0 & C
\end{array}\!\!\right] - r(A) - r(C) \nb
\\
&= r\!\left[\!\! \begin{array}{ccc}
 0  &  I_n & I_n
\\
-B  & A & 0
\\
D  & 0 & C
\end{array}\!\!\right] - r(A) - r(C) \nb
\\
&= r\!\left[\!\! \begin{array}{ccc}
 0  &  I_n & 0
\\
-B  & 0 & -A
\\
D  & 0 & C
\end{array}\!\!\right] - r(A) - r(C) \nb
\\
& =  n + r(N)- r(A) - r(C),
\label{311}
\\
r\!\left[\!\! \begin{array}{cc}
 X_0 - Y_0   &  F_A
\\
F_C & 0
\end{array}\!\!\right] & = r\!\left[\!\! \begin{array}{ccc}
X_0 - Y_0   &  I_n & 0
\\
I_n & 0 & C^*
\\
0 & A & 0
\end{array}\!\!\right]  - r(A) - r(C)  \nb
\\
& = r\!\left[\!\! \begin{array}{ccc}
0  &  I_n & D^*
\\
I_n & 0 & C^*
\\
-B & A & 0
\end{array}\!\!\right]  - r(A) - r(C) \nb
\\
& = r\!\left[\!\! \begin{array}{ccc}
0  &  I_n & 0
\\
I_n & 0 & 0
\\
0 & 0 & BC^* - AD^*
\end{array}\!\!\right]  - r(A) - r(C)  \nb
\\
& = 2n + r(\,  BC^* - AD^* \,) -  r(A) - r(C),
\label{312}
\end{align}
\begin{align}
r\!\left[\!\! \begin{array}{cccc}
 X_0 - Y_0   &  F_A & F_C
\\
F_A & 0 & 0
\end{array}\! \!\right] & = r\!\left[\!\! \begin{array}{cccc}
X_0 - Y_0   &  I_n & I_n & 0
\\
I_n & 0 &  0 &  A^*
\\
0 & A & 0 & 0
\\
0 & 0 & C & 0
\end{array}\!\!\right]  - 2r(A) - r(C)  \nb
\\
& = r\!\left[\!\! \begin{array}{cccc}
0  &  I_n & I_n & 0
\\
I_n & 0 &  0 &  A^*
\\
-B & A & 0 & 0
\\
D & 0 & C & 0
\end{array}\!\!\right]  - 2r(A) - r(C)  \nb
\\
& = r\!\left[\!\! \begin{array}{cccc}
0  &  I_n & 0 & 0
\\
I_n & 0 &  0 &  0
\\
0 & 0 & -A & BA^*
\\
0 & 0 & C & -DA^*
\end{array}\!\!\right]  - 2r(A) - r(C)  \nb
\\
& = 2n + r\!\left[\!\! \begin{array}{cc}
A & BA^*
\\
 C & DA^*
\end{array}\!\!\right]  - 2r(A) - r(C),
\label{313}
\\
r\!\left[\!\! \begin{array}{cccc}
 X_0 - Y_0   &  F_A & F_C
\\
F_C & 0 & 0
\end{array}\! \!\right] & = 2n + r\!\left[\!\! \begin{array}{cc}
A & BC^*
\\
 C & DC^*
\end{array}\!\!\right]  - r(A) - 2r(C),
\label{314}
\\
i_{\pm}\!\left[\!\! \begin{array}{cccc}
 X_0 - Y_0   &  F_A & F_C
\\
F_A & 0 & 0
\\
F_C & 0 & 0
\end{array}\! \!\right] & =
i_{\pm}\!\left[\!\! \begin{array}{ccccc}
  X_0 - Y_0  &  I_n & I_n & 0 & 0
\\
I_n & 0 &  0 &  A^* & 0
\\
I_n & 0 &  0 &  0 & C^*
\\
0 & A & 0 & 0 & 0
\\
0 & 0 & C & 0 & 0
\end{array}\!\!\right]   - r(A) - r(C)  \nb
\\
& =i_{\pm}\!\left[\!\! \begin{array}{ccccc}
0  &  I_n & I_n & -B^*/2 & D^*/2
\\
I_n & 0 &  0 &  A^* & 0
\\
I_n & 0 &  0 &  0 & C^*
\\
-B/2  & A & 0 & 0 & 0
\\
D/2  & 0 & C & 0 & 0
\end{array}\!\!\right]  - r(A) - r(C)  \nb
\\
& =i_{\pm}\!\left[\!\! \begin{array}{ccccc}
0  &  I_n & 0 & 0 & 0
\\
I_n & 0 &  0 &  0 & 0
\\
0 & 0 &  0 &  -A^* & C^*
\\
0  & 0 & -A & AB^* & -AD^*/2
\\
0   & 0 & C & -DA^*/2  & 0
\end{array}\!\!\right]  - r(A) - r(C)  \nb
\\
& = n + i_{\pm}\!\left[\!\! \begin{array}{ccccc}
  0 &  -A^* & C^*
\\
 -A & AB^* & -AD^*/2
\\
C & -DA^*/2  & 0
\end{array}\!\!\right]  - r(A) - r(C)  \nb
\\
& = n + i_{\pm}\!\left[\!\! \begin{array}{ccccc}
 0 &  A^* & C^*
\\
 A & AB^* & 0
\\
 C & 0  & -CD^*
\end{array}\!\!\right]  - r(A) - r(C)  \nb
\\
& = n + i_{\pm}(M)- r(A) - r(C).
\label{315}
\end{align}
Substituting (\ref{311})--(\ref{315}) into (\ref{37})--(\ref{310}) and simplifying leads to (\ref{32})--(\ref{35}).
Results (a)--(g) follow from applying Lemma \ref{T16} to   (\ref{32})--(\ref{35}).
\qquad $\Box$

\medskip

A direct consequence for $AX = B$ and its perturbation equation
$(\,A + \delta A \,)Y = (\,B + \delta B \,)$  is given below.

\begin{corollary} \label{T32}
Assume that both $AX = B$ and its perturbation equation $(\,A + \delta A \,)Y = (\,B + \delta B \,)$
 have Hermitian solutions$,$ and let
\begin{align*}
{\mathcal S} = \{\, X\in {\mathbb C}^{n}_{{\rm H}} \ | \  AX = B \, \}, \ \
{\mathcal T} = \{ \, Y \in {\mathbb C}^{n}_{{\rm H}} \ | \  (\,A + \delta A \,)Y = (\,B + \delta B \,) \, \}.
\end{align*}
Also denote
$$
M = \left[\!\! \begin{array}{ccccc}
 AB^* & 0 & A
 \\
 0 & -(\,A + \delta A \,)(\,B + \delta B \,)^* & (\,A + \delta A \,)
 \\
A^* &  (\,A + \delta A \,)^* & 0
\end{array}\!\!\right], \ \  N = \left[\!\! \begin{array}{cc}
A  & B
\\
\delta A   & \delta B
\end{array}\!\!\right]\!.
$$
Then$,$ the following hold$.$
\begin{enumerate}
\item[{\rm(a)}] There exist $X \in {\mathcal S}$ and $Y\in {\mathcal T}$ such that $X - Y$ is nonsingular
   if and only if
$$
 r\!\left[\!\! \begin{array}{ccc}
A  & B
\\
 \delta A   &  \delta B
\end{array}\!\!\right] =  r(A) + r(\,A + \delta A \,).
$$

\item[{\rm(b)}]  $X - Y$ is nonsingular for all $X \in {\mathcal S}$ and $Y\in {\mathcal T}$ if and only if
$$
2r(N) + r[\, A(\delta B)^* - B(\delta A)^* \,] = r\!\left[\!\! \begin{array}{cc}
A  & BA^*
\\
\delta A  & (\delta B) A^*
\end{array}\!\!\right] +  r\!\left[\!\! \begin{array}{cc}
A  & BA^* + B(\delta A)^*
\\
\delta A  & (\delta B)A^* + (\delta B)(\delta A)^*
\end{array}\!\!\right] +n.
$$

\item[{\rm(c)}] ${\mathcal S} \cap {\mathcal T} \neq \emptyset$ if and only if
$\R\!\left[\!\! \begin{array}{c}
B
\\
\delta B
\end{array}\!\!\right] \subseteq \R\!\left[\!\! \begin{array}{ccc}
A
\\
\delta A
\end{array}\!\!\right]$ and $A(\delta B)^* = B(\delta A)^*.$

\item[{\rm(d)}] There exist $X \in {\mathcal S}$ and $Y\in {\mathcal T}$ such that $X \succ Y$ $(X \prec Y)$
   if and only if $i_{+}(M) = r(A) + r(\,A + \delta A \,)$  $\left( \, i_{-}(M) = r(A) + r(\,A + \delta A \,)  \, \right).$

\item[{\rm(e)}]  $X \succ Y$ $(X \prec Y)$ holds for all $X \in {\mathcal S}$ and $Y\in {\mathcal T}$ if and only if
$i_{-}(M) = r(N) - n$  $\left(\, i_{+}(M) = r(N) - n \, \right)\!.$

\item[{\rm(f)}] There exist $X \in {\mathcal S}$ and $Y\in {\mathcal T}$ such that $X \succcurlyeq Y$
 $(X \preccurlyeq Y)$  if and only if $i_{+}(M) = r(N)$  $\left( \, i_{-}(M) = r(N) \, \right)\!.$

\item[{\rm(g)}]  $X \succcurlyeq Y$ $(X \preccurlyeq Y)$  holds for all $X \in {\mathcal S}$ and $Y\in {\mathcal T}$ if and only if
$$
i_{-}(M) = r(A)  + r(\,A + \delta A \,) -n \ \
\left(\, i_{+}(M)  = r(A)  + r(\,A + \delta A \,) -n  \, \right).
$$
\end{enumerate}
\end{corollary}

\section{Equalities and inequalities for Hermitian solutions of $AXA^* = B$ and its transformed equations}
\renewcommand{\theequation}{\thesection.\arabic{equation}}
\setcounter{section}{4}
\setcounter{equation}{0}

The transformation matrix equation in  (\ref{18}) may reasonably occur in the investigation of (\ref{12}) and its variations.
 For instance,
 \begin{enumerate}
\item[{\rm (a)}] setting $T = A^*$ in (\ref{18}) yields $A^*AXA^{*}A= A^*BA$,
which is the well-known normal equation of (\ref{12}) corresponding to
$\text{trace}[(\, B -  AXA^{*} \,)(\, B -  AXA^{*} \,)^*\,] = \min$;

\item[{\rm (b)}] partitioning $A$ and $B$ in  (\ref{12}) as $A =
\left[\!\! \begin{array}{c} A_1\\ A_2 \end{array} \!\!\right]$ and
 $B = \left[\!\! \begin{array}{cc} B_1 & B_2 \\  B^*_2 & B_3 \end{array} \!\!\right]$,
 and setting  $T =[\, I_{m_1}, \, 0 \,]$ and  $T =[\, 0, \, I_{m_2} \,]$  in (\ref{18}),
 respectively, we obtain two small equations  $A_1XA_1^{*} = B_1$ and $A_2XA_2^{*} = B_3$;

\item[{\rm (c)}] partitioning $A$ and  $X$  in  (\ref{12}) as $A =[\, A_1, \,  A_2 \,]$ and
$X = \left[\!\! \begin{array}{cc} X_1 & X_2 \\  X^*_2 & X_3 \end{array} \!\!\right]$,
 and setting $T =E_{A_1}$ and $T =E_{A_2}$, respectively, we obtain two small equations
$E_{A_1}A_2X_2A^{*}_2E_{A_1} = E_{A_1}BE_{A_1}$ and $E_{A_2}A_1X_1A^{*}_1E_{A_2} = E_{A_2}BE_{A_2}$,
respectively;

\item[{\rm (d)}] decomposing $A$ in (\ref{12}) as a sum  $A = A_1 + A_2$
and setting $T =E_{A_1}$ and $T =E_{A_2}$, respectively,
we obtain two transformed equations $E_{A_1}A_2XA^{*}_2E_{A_1} = E_{A_1}BE_{A_1}$ and
$E_{A_2}A_1XA^{*}_1E_{A_2} = E_{A_2}BE_{A_2}$, respectively.
\end{enumerate}
Since solutions of (\ref{12}) and its transformed
equations are not necessarily the same, it is necessary to consider
relations between the solutions of  (\ref{12}) and its transformed equations.

\begin{theorem} \label{T41}
Assume that {\rm (\ref{12})} has a Hermitian solution$,$  and let
${\cal S}$ and ${\cal T}$ be defined as in {\rm (\ref{19})} and  {\rm (\ref{110})}$.$
Then$,$ the following hold$.$
\begin{enumerate}
\item[{\rm(a)}] ${\cal S} \subseteq  {\cal T}$ always holds$,$ namely$,$ all Hermitian solutions
of {\rm (\ref{12})} are solutions of {\rm (\ref{18})}$.$

\item[{\rm(b)}] ${\cal S} = {\cal T}$ if and only if $r(TA) = r(A),$ or equivalently$,$
$\R(A^*T^*) = \R(A^*).$
 \end{enumerate}
\end{theorem}

\noindent {\bf Proof} \ If $X$ is a Hermitian solution of (\ref{12}), then $TAXA^*T^* = TBT^*$ holds as well.
Hence, we have (a). Note from (\ref{117}) that the set inclusion ${\cal S} \supseteq {\cal T}$ is equivalent
 to the following max-min rank problem
\begin{align}
\max_{Y\in {\cal T}}\,\min_{X\in {\cal S}}r(\, X - Y \,) = 0.
\label{41}
\end{align}
Then, applying (\ref{165}) and simplifying by (\ref{129}), we obtain
\begin{align}
\min_{X\in {\cal S}}\!r(\, X - Y \,) & = \min_{V} r\left[\, A^{\dag}B(A^*)^{\dag}
+ F_{A}V + V^{*}F_{A} - Y \,\right] \nb
\\
& = r\!\left[\!\! \begin{array}{cccc} A^{\dag}B(A^{\dag})^{*} - Y & F_{A} \\ F_{A} & 0
\end{array}\!\!\right] - 2r(F_A)\nb
\\
& = r[\,A^{\dag}B(A^{\dag})^{*} - A^{\dag}AYA^{\dag}A\,] = r(\, B - AYA^*\,).
\label{42}
\end{align}
From Lemma \ref{T111}, the general Hermitian solution of (\ref{18}) can be written as
\begin{align}
 Y = (TA)^{\dag}TBT^*(A^*T^*)^{\dag} - F_{TA}U - U^{*}F_{TA},
\label{43}
\end{align}
where $U$ is arbitrary. Applying (\ref{160}), we obtain
\begin{align}
& \max_{Y\in {\cal T}}\! r(\, B - AYA^*\,) \nb
\\
& = \max_{U} r[\, B - A(TA)^{\dag}TBT^*(A^*T^*)^{\dag}A^* + AF_{TA}UA^* + AU^{*}F_{TA}A^* \,] \nb
\\
& =\min \left\{ r[\, B - A(TA)^{\dag}TBT^*(A^*T^*)^{\dag}A^*,\, A \,],  \ \ r\!\left[\!\!\begin{array}{cc}  B - A(TA)^{\dag}TBT^*(A^*T^*)^{\dag}A^* & AF_{TA}
\\
 F_{TA}A^{*} & 0 \end{array}\!\!\right] \right\} \nb
 \\
& =\min \left\{ r(A),  \ \ r\!\left[\!\!\begin{array}{cc}  B - A(TA)^{\dag}TBT^*(A^*T^*)^{\dag}A^* & AF_{TA}
\\
 F_{TA}A^{*} & 0 \end{array}\!\!\right] \right\}\!.
 \label{44}
 \end{align}
Applying  (\ref{135}) and simplifying by (\ref{120}) and (\ref{123}), we obtain
\begin{align}
& i_{\pm}\!\left[\!\!\begin{array}{cc}  B - A(TA)^{\dag}TBT^*(A^*T^*)^{\dag}A^* & AF_{TA} \\
 F_{TA}A^{*} & 0 \end{array}\!\!\right]  \nb
\\
& = i_{\pm}\!\left[\!\!\begin{array}{ccc}  B - A(TA)^{\dag}TBT^*(A^*T^*)^{\dag}A^* & A & 0
\\
 A^* & 0  & A^*T^*
 \\
 0 & TA & 0 \end{array}\!\!\right] - r(TA) \nb
 \\
 & = i_{\pm}\!\left[\!\!\begin{array}{ccc} B - A(TA)^{\dag}TBT^*(A^*T^*)^{\dag}A^*  & A &
 - BT^* + A(TA)^{\dag}TBT^*
\\
 A^* & 0  & 0
 \\
- TB + TBT^*(A^*T^*)^{\dag}A^* & 0 & 0 \end{array}\!\!\right] - r(TA) \nb
 \\
 & = i_{\pm}\!\left[\!\!\begin{array}{ccc}  0  & A & 0
\\
 A^* & 0  & 0
 \\
0  & 0 & 0\end{array}\!\!\right] - r(TA)  = r(A) - r(TA),
  \label{45}
\end{align}
so that
\begin{align}
& r\!\left[\!\!\begin{array}{cc}  B - A(TA)^{\dag}TBT^*(A^*T^*)^{\dag}A^* & AF_{TA} \\
 F_{TA}A^{*} & 0 \end{array}\!\!\right] = 2r(A) - 2r(TA).
  \label{46}
\end{align}
Substituting (\ref{46}) into (\ref{44}), and then  (\ref{44}) into
 (\ref{42}) yields
\begin{align}
\max_{Y\in {\cal T}}\min_{X\in {\cal S}_1} r(\, X - Y\,) =
\min\{ r(A), \ \ 2r(A) -  2r(TA) \}.
  \label{47}
\end{align}
Finally, substituting (\ref{47}) into (\ref{41}) leads to the result in (b). \qquad $\Box$

\begin{theorem} \label{T42}
Assume that {\rm (\ref{12})} has a  Hermitian solution$,$  and let
${\cal S}$ and  ${\cal T}$  be as given in {\rm (\ref{19})} and {\rm (\ref{110})}$.$
Then$,$  the following hold$.$
\begin{enumerate}
\item[{\rm(a)}] There exist $X \in {\cal S}$ and $Y \in {\cal T}$ such that
$X \succ  Y$ $(X \prec Y)$ if and only if $TA =0.$ So that if  {\rm (\ref{18})} is not null equation$,$ there don't
exist  $X \in {\cal S}$ and $Y \in {\cal T}$ such that
$X \succ  Y$ $(X \prec Y).$

\item[{\rm(b)}] There always exist $X \in {\cal S}$ and $Y \in {\cal T}$ such that
$X \succcurlyeq Y$ $(X \preccurlyeq Y).$
\end{enumerate}
\end{theorem}

\noindent {\bf Proof} \ From  (\ref{118}) and (\ref{119}), there exist $X \in {\cal S}$ and $Y \in {\cal T}$ such that
$X \succ  Y$ $(X \prec Y)$ if and only if
\begin{align}
\max_{X\in {\cal S}, \, Y\in {\cal T}} i_{+}(\, X - Y \,) = n \ \
\left(\max_{X\in {\cal S}, \, Y\in {\cal T}} i_{-}(\, X - Y \,) = n \right);
\label{48}
\end{align}
there exist $X \in {\cal S}$ and $Y \in {\cal T}$ such that $X \succcurlyeq Y$ $(X \preccurlyeq Y)$
if and only if
\begin{align}
\min_{X\in {\cal S}, \, Y\in {\cal T}} i_{-}(\, X - Y \,) = 0 \ \
\left(\min_{X\in {\cal S}_1, \, Y\in {\cal T}} i_{+}(\, X - Y \,) = 0 \right).
\label{49}
\end{align}
From (\ref{139}) and (\ref{43}), the difference of $X - Y$ can be written as
\begin{align}
X -  Y & =  A^{\dag}B(A^*)^{\dag} -(TA)^{\dag}TBT^*(A^*T^*)^{\dag} + F_{A}V + V^{*}F_{A}
+ F_{TA}U + U^{*}F_{TA} \nb
\\
& =  A^{\dag}B(A^*)^{\dag} -(TA)^{\dag}TBT^*(A^*T^*)^{\dag} + [\, F_{A}, \, F_{TA}\,]\!\left[\!\!\begin{array}{c}
 V \\ U \end{array}\!\!\right] + [\, V^*, \, U^*\,]\!\left[\!\!\begin{array}{c}
 F_A \\ F_{TA} \end{array}\!\!\right]\!.
\label{410}
\end{align}
Applying (\ref{166}) and simplifying by (\ref{128}) and (\ref{45}), we obtain
\begin{align}
\max_{X\in {\cal S}, \, Y\in {\cal T}}\!i_{\pm}(\, X - Y \,)
& = i_{\pm}\!\left[\!\! \begin{array}{cccc} A^{\dag}B(A^*)^{\dag} -(TA)^{\dag}TBT^*(A^*T^*)^{\dag} & F_{A}
&  F_{TA} \\ F_{A} & 0 & 0 \\  F_{TA} & 0 & 0
\end{array}\! \!\right] \nb
\\
& = i_{\pm}\!\left[\!\! \begin{array}{cccc} B - A(TA)^{\dag}TBT^*(A^*T^*)^{\dag}A^* &  AF_{TA}
 \\ F_{A}A^* & 0 \end{array}\! \!\right] + r(F_A)  \nb
 \\
& = r(A) - r(TA)  + n - r(A) =  n - r(TA).
\label{411}
\end{align}
Substituting (\ref{411}) into (\ref{48}) leads to the result in (a).

Applying (\ref{167}) and simplifying by (\ref{125}), (\ref{128}) and (\ref{45}), we obtain
\begin{align}
\min_{X\in {\cal S}, \, Y\in {\cal T}}\!i_{\pm}(\, X - Y \,)
& = i_{\pm}\!\left[\!\! \begin{array}{cccc} A^{\dag}B(A^*)^{\dag} -(TA)^{\dag}TBT^*(A^*T^*)^{\dag} & F_{A}
&  F_{TA} \\ F_{A} & 0 & 0 \\  F_{TA} & 0 & 0
\end{array}\! \!\right] - r[\, F_{A}, \,  F_{TA} \,] \nb
\\
& = i_{\pm}\!\left[\!\! \begin{array}{cccc} B - A(TA)^{\dag}TBT^*(A^*T^*)^{\dag}A^* &  AF_{TA}
 \\ F_{A}A^* & 0 \end{array}\! \!\right]  - r(AF_{TA}) \nb
 \\
& = r(A) - r(TA) - r(A) + r(TA) = 0.
\label{412}
\end{align}
Substituting (\ref{412}) into (\ref{49}) leads to the result in (b). \qquad $\Box$

\section{Average equalities for Hermitian solutions of $AXA^{*} = B$ and its two transformed equations}
\renewcommand{\theequation}{\thesection.\arabic{equation}}
\setcounter{section}{5}
\setcounter{equation}{0}

In this section, we study the relations between the two sets in (\ref{111}) and (\ref{112}).

\begin{theorem} \label{T51}
Assume that {\rm (\ref{12})} has a Hermitian solution$,$  and let
${\cal S}$ and ${\cal T}$ be defined as in {\rm (\ref{111})} and {\rm (\ref{112})}$,$ in which
 $T_1 \in \mathbb C^{p_1\times n}$ and $T_2 \in \mathbb C^{p_2\times n}.$
Then$,$
\begin{enumerate}
\item[{\rm(a)}] ${\cal S} \subseteq  {\cal T}$ always holds$.$

\item[{\rm(b)}] ${\cal S} = {\cal T}$ if and only if
\begin{align}
r\!\left[\!\!\begin{array}{cc}  T_1A
\\
T_2A
\end{array}\!\!\right] =  r(T_1A) + r(T_2A) - r(A).
\label{51}
\end{align}

\item[{\rm(c)}] In particular$,$ if $r(T_1A) = r(T_2A) = r(A),$ then ${\cal S} = {\cal T}.$

 \end{enumerate}
\end{theorem}

\noindent {\bf Proof} \ If $X$ is a Hermitian solution of (\ref{12}), then
$T_1AXA^*T^*_1 = T_1BT^*_1$ and $T_2AXA^*T^*_2 = T_2BT^*_2$ hold as well, and
\begin{align}
r\!\left[\!\!\begin{array}{ccccc}
T_1BT_1^* & 0 &   T_1A
\\
0 & -T_2BT_2^* & T_2A
\\
A^*T^*_1  & A^*T^*_2 & 0\end{array}\!\!\right] =  2r\!\left[\!\!\begin{array}{ccccc}
T_1A
\\
T_2A
\end{array}\!\!\right]\!.
\label{52}
\end{align}
holds by Lemma \ref{T112}(c). This fact means that any Hermitian solution $X$ of (\ref{12}) can be written
as $X = (\, X + X \,)/2$, the average of the two Hermitian solutions of the two equations in (\ref{112}).
 Hence, we have (a).

 Note from  (\ref{117}) that the set inclusion ${\cal S} \supseteq {\cal T}$ is equivalent to
\begin{align}
\max_{Y\in {\cal T}}\,\min_{X\in {\cal S}}r(\, X - Y \,) = 0.
\label{53}
\end{align}
Applying (\ref{165}), we first obtain
\begin{align}
\min_{X\in {\cal S}}\!r(\, X - Y \,)  = \min_{V} r\left[\, A^{\dag}B(A^*)^{\dag}
+ F_{A}V + V^{*}F_{A} - Y \,\right]  = r(\, B - AYA^*\,).
\label{54}
\end{align}
From Lemma \ref{T111}, the general expression of the matrices of the two equations in (\ref{112})
 can be written as
\begin{align}
 Y & = (T_1A)^{\dag}T_1BT^*_1(A^*T^*_1)^{\dag}/2  + (T_2A)^{\dag}T_2BT^*_2(A^*T^*_2)^{\dag}/2 \nb
\\
& \ \  - F_{T_1A}U_1 - U^{*}_1F_{T_1A}  - F_{T_2A}U_2 - U^{*}_2F_{T_2A},
\label{55}
\end{align}
where $U_1$ and $U_2$ are arbitrary. Then
\begin{align}
 B - AYA^* & = B - A(T_1A)^{\dag}T_1BT^*_1(A^*T^*_1)^{\dag}A^*/2  + A(T_2A)^{\dag}T_2BT^*_2(A^*T^*_2)^{\dag}A^*/2 \nb
  \\
  &  \ \ \ - AF_{T_1A}U_1A^* - AU^{*}_1F_{T_1A}A^*  - AF_{T_2A}U_2A^* - AU^{*}_2F_{T_2A}A^* \nb
\\
& = G  - [\, AF_{T_1A}, \,  AF_{T_2A} \,]\!\left[\!\begin{array}{c}
 U_1 \\ U_2 \end{array} \!\right]A^*  - A[\, U_1^*, \,  U_2^*\,]\left[\!\begin{array}{c}
 F_{T_1A}A^* \\  F_{T_2A}A^*  \end{array}\!\right]\!,
\label{56}
\end{align}
where $ G = B - A(T_1A)^{\dag}T_1BT^*_1(A^*T^*_1)^{\dag}A^*/2 - A(T_2A)^{\dag}T_2BT^*_2(A^*T^*_2)^{\dag}A^*/2$.
Applying (\ref{160}) gives
\begin{align}
 \max_{Y\in {\cal T}} r(\, B - AYA^*\,)  & = \max_{U_1,\, U_2} r\!\left( G  -
 [\, AF_{T_1A}, \,  AF_{T_2A} \,]\!\left[\!\!\begin{array}{c}
 U_1 \\ U_2 \end{array}\!\!\right]A^*  - A[\, U_1^*, \,  U_2^*\,]\!\left[\!\!\begin{array}{c}
 F_{T_1A}A^* \\  F_{T_2A}A^*  \end{array}\!\!\right] \right) \nb
\\
& =\min \!\left\{ r(A), \ \ r\!\left[\!\!\begin{array}{ccc}  G & AF_{T_1A} & AF_{T_2A}
\\
 F_{T_1A}A^* & 0 & 0
 \\
 F_{T_2A}A^* & 0 & 0
  \end{array}\!\!\right] \right\}\!,
 \label{57}
 \end{align}
 where applying (\ref{135}) and simplifying by elementary matrix operations, we obtain
\begin{align}
& r\!\left[\!\!\begin{array}{ccc}  G & AF_{T_1A} & AF_{T_2A}
\\
 F_{T_1A}A^* & 0 & 0
 \\
 F_{T_2A}A^* & 0 & 0
  \end{array}\!\!\right] \nb
  \\
& = r\!\left[\!\!\begin{array}{ccccc}
2B - A(T_1A)^{\dag}T_1BT^*_1(A^*T^*_1)^{\dag}A^* - A(T_2A)^{\dag}T_2BT^*_2(A^*T^*_2)^{\dag}A^* &\! A &\! A &\! 0 & \! 0
\\
 A^* & \!0 & \!0  & \!A^*T^*_1 & \! 0
 \\
 A^* & \!0 & \!0 & \!0 & \!A^*T^*_2
\\
 0 & \!T_1A & \!0 & \!0 & \!0
 \\
0 & \!0 & \!T_2A & \!0 & \! 0
\end{array}\!\!\right] \nb
\\
& \ \ \  - 2r(T_1A) - 2r(T_2A) \nb
\\
& = r\!\left[\!\!\begin{array}{ccccc}
0 & A & 0 & 0 & 0
\\
 A^* & 0 & 0  & A^*T^*_1 & 0
 \\
 A^* & 0 & 0 & 0 & A^*T^*_2
\\
 -2T_1B + T_1BT^*_1(A^*T^*_1)^{\dag}A^*  + T_1A(T_2A)^{\dag}T_2BT^*_2(A^*T^*_2)^{\dag}A^* & 0 & -T_1A  & 0 & 0
 \\
0 & 0 &  T_2A  & 0 & 0
\end{array}\!\!\right]  \nb
\\
&  \ \ \ \  - 2r(T_1A) - 2r(T_2A)  \nb
\\
& = r\!\left[\!\!\begin{array}{ccccc}
 A^* & 0   & 0 & 0
 \\
0 & 0 &  -A^*T^*_1  & A^*T^*_2
\\
 0 & -T_1A  & 2T_1BT_1^* - T_1BT^*_1  -T_1A(T_2A)^{\dag}T_2BT^*_2(A^*T^*_2)^{\dag}A^*T_1^* & 0
 \\
0 &   T_2A  & 0 & 0
\end{array}\!\!\right]  \nb
\\
&   \ \ \ \ \  + r(A)  - 2r(T_1A) - 2r(T_2A)  \nb
\\
& = r\!\left[\!\!\begin{array}{ccccc}
 0 &  -A^*T^*_1  & A^*T^*_2
\\
  -T_1A  & T_1BT_1^*   -T_1A(T_2A)^{\dag}T_2BT^*_2(A^*T^*_2)^{\dag}A^*T_1^* & 0
 \\
   T_2A  & 0 & 0
\end{array}\!\!\right]  + 2r(A)  - 2r(T_1A) - 2r(T_2A)  \nb
\\
& = r\!\left[\!\!\begin{array}{ccccc}
 0 &  -A^*T^*_1  & A^*T^*_2
\\
  -T_1A  & T_1BT_1^*  & -T_1A(T_2A)^{\dag}T_2BT^*_2
 \\
   T_2A  & 0 & 0
\end{array}\!\!\right]  + 2r(A)  - 2r(T_1A) - 2r(T_2A)  \nb
\\
& = r\!\left[\!\!\begin{array}{ccccc}
 0 &  -A^*T^*_1  & A^*T^*_2
\\
  -T_1A  & T_1BT_1^*  & 0
   \\
   T_2A  & 0 & -T_2BT_2^*
\end{array}\!\!\right]  + 2r(A)  - 2r(T_1A) - 2r(T_2A)  \nb
\\
& = 2r\!\left[\!\!\begin{array}{cc}
 T_1A
\\
T_2A
\end{array}\!\!\right]  + 2r(A)  - 2r(T_1A) - 2r(T_2A) \ \ (\mbox{by (\ref{52})}).
\label{58}
\end{align}
Hence,
\begin{align}
 \max_{Y\in {\cal T}}\! r(\, B - AYA^*\,)  =\min\! \left\{\! r(A),  \ 2r\!\left[\!\!\begin{array}{cc}
 T_1A
\\
T_2A
\end{array}\!\!\right] \! + 2r(A)  - 2r(T_1A) - 2r(T_2A) \!\right\}\!.
 \label{59}
 \end{align}
Setting the both sides to zero leads to the equality in (\ref{51}).  Under
$r(T_1A) = r(T_2A) = r(A)$, both sides of  (\ref{51}) reduces to $2r(A)$. Hence,
(c) follows.  \qquad $\Box$

\medskip

We next give some consequences of Theorem \ref{T51}. Partitioning $A$ and $B$ in  (\ref{12}) as
$$
\left[\!\! \begin{array}{c} A_1\\ A_2 \end{array} \!\!\right]X[\,  A^*_1, \, A^*_2 \,] = \left[\!\! \begin{array}{cc} B_1 & B_3 \\  B^*_3 & B_2 \end{array} \!\!\right]\!,
$$
where $A_1\in \mathbb C^{m_1\times n}$ $A_2\in \mathbb C^{m_2\times n}$, $B_1 \in {\mathbb C}_{{\rm H}}^{m_1}$,
  $B_2 \in {\mathbb C}_{{\rm H}}^{m_2},$  $B_3\in \mathbb C^{m_1\times m_2}$, $m_1 + m_2 = m,$ and setting
  $T_1 =[\, I_{m_1}, \, 0 \,]$ and  $T =[\, 0, \, I_{m_2} \,]$  in (\ref{18}), respectively,
 we obtain two small equations
$$
A_1XA_1^{*} = B_1 \ \  {\rm and} \ \ A_2XA_2^{*} = B_2.
$$
 Let
\begin{align}
 {\cal T} = \left\{\,  (X_1 + X_2)/2 \ | \ A_1X_1A_1^{*} = B_1, \
  A_2X_2A_2^{*} = B_2, \ X_1, \, X_2 \in {\mathbb C}_{{\rm H}}^{n} \,  \right\}.
\label{510}
\end{align}
Applying Theorem \ref{T51} to  (\ref{111}) and (\ref{510}), we obtain the following result.

\begin{corollary} \label{T53}
Assume that {\rm (\ref{12})} has a solution$,$  and let ${\cal S}$ be as given in
{\rm (\ref{111})}  and  ${\cal T}$ as in {\rm (\ref{510})}$.$ Then$,$ ${\cal S} = {\cal T}$ if and only if
$\R(A^*_1) = \R(A^*_2).$
\end{corollary}

Decomposing $A$ in (\ref{12}) as $A = A_1 + A_2$ and setting $T =E_{A_1}$ and $T =E_{A_2}$
respectively  yields the following two transformed equations
$$
E_{A_1}A_2XA^{*}_2E_{A_1} = E_{A_1}BE_{A_1}, \ \ E_{A_2}A_1XA^{*}_1E_{A_2} = E_{A_2}BE_{A_2}.
$$
Also let
\begin{align}
 {\cal T} = \left\{(\,X_1 + X_2\,)/2  \,|
  \, E_{A_2}\!A_1X_1A^{*}_1E_{A_2} = E_{A_2}\!BE_{A_2},  \,
 E_{A_1}\!A_2X_2A^{*}_2E_{A_1} = E_{A_1}\!BE_{A_1}, \ X_1, \, X_2 \in {\mathbb C}_{{\rm H}}^{n} \right\}.
\label{511}
\end{align}

\begin{corollary} \label{T54}
Assume that  {\rm (\ref{12})} has a solution$,$ and let
${\cal S}$ and ${\cal T}$ be as given in {\rm (\ref{111})} and {\rm (\ref{511})}$.$
  Then$,$ the following hold$.$
  \begin{enumerate}
\item[{\rm(a)}] ${\cal S} = {\cal T}$ if and only if
\begin{align}
r\!\left[\!\! \begin{array}{cccc}
A_1 & 0 & A_2
\\
0 & A_2& A_1
\end{array}\! \!\right] = 2r[\,A_1,\,A_2\,] - r(A).
\label{512}
\end{align}

\item[{\rm(b)}]  Under ${\mathscr R}(A_1) \cap {\mathscr R}(A_2) =\{ 0\},$
 ${\cal S} = {\cal T}$  if and only if ${\mathscr R}(A_1^*) ={\mathscr R}(A_2^*).$
 \end{enumerate}

  \end{corollary}

\noindent {\bf Proof} \ From Theorem \ref{T51}(b), ${\cal S} = {\cal T}$ if and only if
\begin{align}
r\!\left[\!\!\begin{array}{ccccc}
 E_{A_2}A_1
\\
 E_{A_1}A_2
\end{array}\!\!\right] =  r(E_{A_2}A_1) + r(E_{A_1}A_2) - r(A),
\label{513}
\end{align}
which is equivalent to  (\ref{512}) by (\ref{124}). \qquad $\Box$

\medskip

In addition to the average equalities of Hermitian solutions of (\ref{12}) and its two transformed equations,
it would be of interest to consider weighted average equalities for the Hermitian solutions of
(\ref{12}) and its $k$ transformed equations
\begin{align}
T_iAX_iA^*T^*_i = T_iBT^*_i,  \ \ i = 1, \ldots, k,
\label{514}
\end{align}
where $T_i \in \mathbb C^{p_i\times n}$.  As usual, define
\begin{align}
{\cal T} & = \left\{\sum_{i =1}^{k}\lambda_iX_i  \,  \left|  \, T_iAX_iA^*T^*_i = T_iBT^*_i,
 \,  X_i = X^*_i, \,  \sum_{i =1}^{k}\lambda_i =1, \ \lambda_i >0, \,   i = 1, \ldots, k  \right. \right\}\!.
\label{515}
\end{align}
 An open problem is to establish necessary and sufficient conditions for  ${\cal S} = {\cal T}$ to hold.

\section{Equalities and inequalities between least-squares and least-rank Hermitian solutions of $AXA^{*} =B$}
\setcounter{section}{6}
\setcounter{equation}{0}

It is well known that the normal equation corresponding to the norm minimization problem in (\ref{113}) is given by
\begin{equation}
A^{*}AXA^{*}A = A^{*}BA,
\label{61}
\end{equation}
see \cite{BG}, while  the normal equation corresponding to the rank minimization problem in {\rm (\ref{114})} is given by
\begin{equation}
 E_{T_1}YE_{T_1} =  - E_{T_1}TM^{\dag}T^{*}E_{T_1},
\label{62}
\end{equation}
where $M = \left[\!\!\begin{array}{cc} B & A \\ A^{*} & 0 \end{array}\!\!\right]\!,$  $T = [\, 0, \, I_n \,]$ and $T_1 = TF_M$;
see \cite{T-laa10}.  Both (\ref{61}) and (\ref{62}) are transformed equations of $AXA^* = B$. Relations between solutions
of the two equations were considered in \cite{T-nla}. In this section, we
 establish a group of formulas for calculating the extremal ranks and inertias of $X - Y$ for their Hermitian solutions,
 and  use the formulas to solve Problem \ref{T15}.

\begin{lemma}\label{T61}
Let  $A\in \mathbb C^{m\times n}$ and $B \in {\mathbb C}_{{\rm H}}^{m}$ be given$.$ Then the following hold$.$
\begin{enumerate}
\item[{\rm(a)}]  {\rm \cite{BG}} The solution of {\rm (\ref{113})} is
\begin{align}
\min_{X \in \mathbb C^{n}_{{\rm H}}}\|\, B - AXA^{*} \,\|_{F} &  =  \|\,B - AA^{\dag}BAA^{\dag}\,\|_{F},
\label{63}
\\
{\rm arg}\!\!\min_{X \in \mathbb C^{n}_{{\rm H}}} \|\, B - AXA^{*} \,\|_{F}  & = A^{\dag}B(A^{\dag})^{*} + F_AU + U^{*}F_A,
\label{64}
\end{align}
where $U \in \mathbb C^{n \times n}$ is arbitrary$.$

\item[{\rm(b)}] {\rm \cite{T-laa10}}  The solution of {\rm (\ref{114})} is
\begin{align}
\min_{Y \in \mathbb C^{n}_{{\rm H}}}\!r( \, B - AYA^{*} \,)  & = 2r[\, A, \, B \, ] -
 r\!\left[\!\!\begin{array}{cc} B & A \\ A^{*} & 0 \end{array}\!\!\right]\!,
\label{65}
\\
{\rm arg}\!\!\min_{Y \in \mathbb C^{n}_{{\rm H}}}\!r( \, B - AYA^{*} \,)  & = - TM^{\dag}T^{*} + T_1V + V^{*}T^{*}_1,
\label{66}
\end{align}
where $V\in \mathbb C^{(m + n)\times n}$ is arbitrary.
\end{enumerate}
\end{lemma}

\begin{theorem} \label{T62}
Let ${\cal S}$ and ${\cal T}$ be as given in {\rm (\ref{113})} and {\rm (\ref{114})}$,$
and define
$$
M = \left[\!\!\begin{array}{ccccc} B &  A
\\
A^{*} &  0 \end{array} \!\!\right], \ \ N = \left[\!\!\begin{array}{ccccc} B & BA & A
\\
A^{*}B &  0  &  0
\\
A^{*} & 0 & 0
 \end{array} \!\!\right]\!.
$$
Then$,$
\begin{align}
\max_{X \in {\mathcal S},  Y \in {\mathcal T}}r(\,X - Y\,) & = \min \{ \, n, \ \ 2n + r(N) - 2r(A) -  r(M) \},
\label{67}
\\
\min_{X \in {\mathcal S},  Y \in {\mathcal T}} r(\,X - Y\,) & = r(N)  + r(M) - 2r[\, A, \, B \,] - 2r(A),
\label{68}
\\
\max_{X \in {\mathcal S},  Y \in {\mathcal T}}i_{\pm}(\,X - Y\,) & = i_{\mp}(N) + n - r(A) -  i_{\mp}(M),
\label{69}
\\
\min_{X \in {\mathcal S},  Y \in {\mathcal T}} i_{\pm}(\,X - Y\,) & = i_{\mp}(N)  + i_{\pm}(M) - r[\, A, \, B \,] - r(A)
\label{610}
\end{align}
hold$.$ Under the condition $B \succcurlyeq 0,$
\begin{align}
\max_{X \in {\mathcal S},  Y \in {\mathcal T}}r(\,X - Y\,) & = \min \{ \, n, \ \ 2n + r[\, A, \, BA \,]  - 3r(A) \},
\label{611}
\\
\min_{X \in {\mathcal S},  Y \in {\mathcal T}} r(\,X - Y\,) & = r[\, A, \, BA \,] - r(A),
\label{612}
\\
\max_{X \in {\mathcal S},  Y \in {\mathcal T}}i_{+}(\,X - Y\,) & = r[\, A, \, BA \,]  + n - 2r(A),
\label{613}
\\
\max_{X \in {\mathcal S},  Y \in {\mathcal T}} i_{-}(\,X - Y\,) & =  n - r(A),
\label{614}
\\
\min_{X \in {\mathcal S},  Y \in {\mathcal T}} i_{+}(\,X - Y\,) & =
r[\, A, \, BA \,]  - r(A),
\label{615}
\\
\min_{X \in {\mathcal S},  Y \in {\mathcal T}} i_{-}(\,X - Y\,) & = 0.
\label{616}
\end{align}
In consequence$,$ the following hold$.$
\begin{enumerate}
\item[{\rm (a)}]  There exist $X \in {\cal S}$ and $Y \in {\cal T}$ such that $X -Y$ is nonsingular
if and only if  $r(N) \geqslant  2r(A) + r(M) -n.$

\item[{\rm (b)}]   $X -Y$ is nonsingular for all $X \in {\cal S}$ and $Y \in {\cal T}$
 if and only if  $r(N)  + r(M) = 2r[\, A, \, B \,] + 2r(A) +n.$

\item[{\rm (c)}] There exist $X \in {\cal S}$ and $Y \in {\cal T}$ such that $X = Y$  if and only if
 $r(N)  + r(M) = 2r[\, A, \, B \,] + 2r(A).$

\item[{\rm (d)}] There exist $X \in {\cal S}$ and $Y \in {\cal T}$ such that $X \succ  Y$  if and only if  $i_{-}(N) =  i_{-}(M) + r(A).$

\item[{\rm (e)}] $X \succ  Y$ holds for all $X \in {\cal S}$ and $Y \in {\cal T}$ if and only if
$i_{-}(N) = r(A) + r[\, A, \, B \,] - i_{+}(M) +n.$

\item[{\rm (f)}] There exist $X \in {\cal S}$ and $Y \in {\cal T}$ such that $X \succcurlyeq Y$ if and only if $i_{+}(N) = r(A) + r[\, A, \,B \,] - i_{-}(M).$

\item[{\rm (g)}] $X \succcurlyeq Y$ holds for all $X \in {\cal S}$ and $Y \in {\cal T}$ if and only if  $i_{+}(N) = i_{+}(M) + r(A) -n.$

\item[{\rm (h)}] There exist $X \in {\cal S}$ and $Y \in {\cal T}$ such that $X \prec Y$ if and only if  $i_{+}(N) =  i_{+}(M) + r(A).$

\item[{\rm (i)}]  $X \prec Y$ holds for all $X \in {\cal S}$ and $Y \in {\cal T}$ if and only if $i_{+}(N) = r(A) + r[\, A, \, B \,] - i_{-}(M) +n.$

\item[{\rm (j)}] There exist $X \in {\cal S}$ and $Y \in {\cal T}$ such that $X \preccurlyeq Y$ if and only if
 $i_{-}(N) = r(A) + r[\, A, \,B \,] - i_{+}(M).$

\item[{\rm (k)}] $X \preccurlyeq Y$ holds for all $X \in {\cal S}$ and $Y \in {\cal T}$ if and only if  $i_{-}(N) = i_{-}(M) + r(A) -n.$
\end{enumerate}
Under the condition $B \succcurlyeq 0,$  the following hold$.$
\begin{enumerate}
\item[{\rm (l)}]  There exist $X \in {\cal S}$ and $Y \in {\cal T}$ such that $X -  Y$ is  nonsingular  if and only if
 $r[\, A, \, BA \,] \geqslant  3r(A) -n.$

\item[{\rm (m)}]   $X -Y$ is nonsingular for all $X \in {\cal S}$ and $Y \in {\cal T}$
 if and only if  $r[\, A, \, BA \,] = r(A) +n.$

\item[{\rm (n)}]  There exist $X \in {\cal S}$ and $Y \in {\cal T}$ such that $X =  Y$ if and only if
 $\R(BA) \subseteq \R(A).$

\item[{\rm (o)}]  There exist $X \in {\cal S}$ and $Y \in {\cal T}$ such that $X \succ  Y$ if and only if
 $r[\, A, \,BA \,] = 2r(A).$

\item[{\rm (p)}]   $X \succ  Y$ holds for all $X \in {\cal S}$ and $Y \in {\cal T}$ if and only if
 $r[\, A, \, BA \,] = r(A) + n.$

\item[{\rm (q)}]  There always exist $X \in {\cal S}$ and $Y \in {\cal T}$ such that $X \succcurlyeq Y.$

\item[{\rm (r)}] $X \succcurlyeq Y$ holds for all $X \in {\cal S}$ and $Y \in {\cal T}$ if and only if
 $ r(A) =n.$

\item[{\rm (s)}]  There  exist $X \in {\cal S}$ and $Y \in {\cal T}$ such that $X \prec Y$ if and only if
 $A = 0.$

\item[{\rm (t)}] There exist $X \in {\cal S}$ and $Y \in {\cal T}$ such that $X \preccurlyeq Y$ if and only if
  $\R(BA) \subseteq \R(A).$

\item[{\rm (u)}] $X \preccurlyeq Y$ holds for all $X \in {\cal S}$ and $Y \in {\cal T}$  if and only if
  $r[\, A, \,BA \,] = 2r(A) -n.$
\end{enumerate}
\end{theorem}

\noindent {\bf Proof} \ Applying  (\ref{142})--(\ref{145})  to (\ref{61}) and (\ref{62}) gives
\begin{align}
\max_{X \in {\mathcal S},  Y \in {\mathcal T}} r(\,X - Y\,) &
 = \left\{ n, \  r\!\left[\!\!\begin{array}{ccc} A^{*}BA & 0 & A^{*}A \\ 0  & E_{T_1}TM^{\dag}T^{*}E_{T_1} & E_{T_1}
\\ A^{*}A & E_{T_1} & 0 \end{array} \!\!\right] + 2n  - 2r(A) - 2r(E_{T_1}) \right\}\!,
\label{617}
\\
\min_{X \in {\mathcal S},  Y \in {\mathcal T}}
r(\,X - Y\,) & = r\!\left[\!\!\begin{array}{ccc} A^{*}BA & 0 & A^{*}A \\ 0
& E_{T_1}TM^{\dag}T^{*}E_{T_1} & E_{T_1}
\\ A^{*}A & E_{T_1} & 0
 \end{array} \!\!\right] - 2r[\, A^{*}A, \,  E_{T_1}\,],
\label{618}
\\
\max_{X \in {\mathcal S},  Y \in {\mathcal T}}i_{\pm}(\,X - Y\,) &
 = i_{\pm}\!\left[\!\!\begin{array}{ccc} A^{*}BA & 0 & A^{*}A \\ 0  & E_{T_1}TM^{\dag}T^{*}E_{T_1} & E_{T_1}
\\ A^{*}A & E_{T_1} & 0 \end{array} \!\!\right]  + n - r(A) - r(E_{T_1}),
\label{619}
\\
\min_{X \in {\mathcal S},  Y \in {\mathcal T}}
i_{\pm}(\,X - Y\,) & = i_{\pm}\!\left[\!\!\begin{array}{ccc} A^{*}BA & 0 & A^{*}A \\ 0
& E_{T_1}TM^{\dag}T^{*}E_{T_1} & E_{T_1}
\\ A^{*}A & E_{T_1} & 0
 \end{array} \!\!\right] - r[\, A^{*}A, \,  E_{T_1}\,].
\label{620}
\end{align}
Applying Lemma \ref{T18} and simplifying by elementary operations of matrices that
\begin{align}
r( E_{T_1}) & = n - r(T_1) = n - r(TF_M) = n -
 r\!\left[\!\!\begin{array}{cc} 0 &  I_n  \\  B  & A \\ A^{*} & 0 \end{array} \!\!\right] + r(M)  \nonumber
 \\
 & = r(M) - r[\, A, \, B \,],
\label{621}
\\
r[\,  A^{*}A, \, E_{T_1} \,] & = r(A^{*}AT_1) + r(E_{T_1}) = r(ATF_M) + r(M) - r[\, A, \, B \,] \nonumber
\\
& = r\!\left[\!\!\begin{array}{c} A^{*}AT \\ M \end{array} \!\!\right]  - r[\, A, \, B \,]
= r\!\left[\!\!\begin{array}{cc} 0&  A  \\  B  & A \\ A^{*} & 0 \end{array} \!\!\right] - r[\, A, \, B \,] = r(A),
\label{622}
\end{align}
and by Lemma \ref{T18} and elementary operations of matrices that
\begin{align}
& i_{\pm}\!\left[\!\!\begin{array}{ccc} A^{*}BA & 0 & A^{*}A \\ 0  & E_{T_1}TM^{\dag}T^{*}E_{T_1} & E_{T_1}
\\ A^{*}A & E_{T_1} & 0
 \end{array} \!\!\right] \nonumber
 \\
&  = i_{\pm}\!\left[\!\!\begin{array}{ccccc} A^{*}BA & 0 & A^{*}A & 0 \\ 0  & TM^{\dag}T^{*} & I_n & TF_M
\\
A^{*}A & I_n & 0 & 0
\\
0 & E_MT & 0 & 0
 \end{array}\!\!\right] - r(TF_M) \ \ \ \mbox{ (by {\rm (\ref{135})})} \nonumber
\\
& = i_{\pm}\!\left[\!\!\begin{array}{ccccc} A^{*}BA + A^{*}ATM^{\dag}T^{*}A^{*}A & 0 & 0 & -A^{*}ATE_M \\ 0  &  0 & I_n & 0
\\
0 & I_n & 0 & 0
\\
-E_MT^{*}A^{*}A & 0 & 0 & 0
 \end{array} \!\!\right]  - r\!\left[\!\!\begin{array}{c} M\\ T \end{array} \!\!\right]+ r(M)   \ \ \ \  \mbox{(by {\rm (\ref{120})})}\nonumber
\\
& = i_{\pm}\!\left[\!\!\begin{array}{ccccc} A^{*}BA + A^{*}ATM^{\dag}T^{*}A^{*}A & A^{*}ATE_M
\\
E_MT^{*}A^{*}A & 0 \end{array}\!\!\right]  - r[\, A, \, B \,]  + r(M) \ \ \ \mbox{ (by {\rm (\ref{122})})}  \nonumber
\\
& =  i_{\pm}\!\left[\!\!\begin{array}{ccccc} A^{*}BA & A^{*}AT
\\ T^{*}A^{*}A & - M
 \end{array}\!\!\right]  - r[\, A, \, B \,] + i_{\pm}(M)  \ \ \ \mbox{(by {\rm (\ref{128})})}  \nonumber
\\
& = i_{\pm}\!\left[\!\!\begin{array}{ccccc} A^{*}BA & 0 & A^{*}A
\\
0 &  -B  &  -A
\\
A^{*}A & -A^{*} & 0
 \end{array} \!\!\right]  -  r[\, A, \, B \,] + i_{\pm}(M)  \nonumber
\\
& = i_{\pm}\!\left[\!\!\begin{array}{ccccc} 0 & -A^{*}B & 0
\\
-BA &  -B  &  -A
\\
0 & -A^{*} & 0
 \end{array} \!\!\right]  - r[\, A, \, B \,] + i_{\pm}(M)   \ \ \ \  \mbox{(by {\rm (\ref{120})})} \nonumber
 \\
 & = i_{\mp}\!\left[\!\!\begin{array}{ccccc} B & BA & A
\\
A^{*}B &  0  &  0
\\
A^{*} & 0 & 0
 \end{array} \!\!\right] - r[\, A, \, B \,] + i_{\pm}(M) \ \ \ \  \mbox{(by {\rm (\ref{120})} and {\rm (\ref{121})})}
 \\
 & = i_{\mp}(N) - r[\, A, \, B \,] + i_{\pm}(M).
\label{623}
\end{align}
Substituting (\ref{621})--(\ref{623}) into (\ref{617})--(\ref{620}) yields (\ref{67})--(\ref{610}).

Applying (\ref{130}) to (\ref{67})--(\ref{610}) yields (\ref{611})--(\ref{617}). Applying Lemma \ref{T16} to
(\ref{67})--(\ref{610}) leads to (a)--(k). Under the condition $B \succcurlyeq 0,$
applying Lemma \ref{T16} to (\ref{611})--(\ref{617}) leads to (l)--(u). \qquad $\Box$

\medskip

The positive semi-definite least-squares Hermitian solution of (\ref{12}) is defined to an $X$ that satisfies
\begin{equation}
 \|\, B - AXA^{*} \,\|_F  = \min \ \ {\rm s.t.} \ \ X \succcurlyeq 0;
\label{624}
\end{equation}
the positive semi-definite least-rank Hermitian solution of (\ref{12}) is defined to be a $Y$ that satisfies
\begin{equation}
 r(\,B - AYA^{*}\,)  = \min \ \ {\rm s.t.} \ \ Y \succcurlyeq 0.
\label{625}
\end{equation}
The normal equations corresponding to  (\ref{624}) and (\ref{625}) are given by
\begin{align}
& A^{*}AXA^{*}A = A^{*}BA,  \ \   X \succcurlyeq 0,
\label{626}
\\
&  E_{T_1}YE_{T_1} =  - E_{T_1}TM^{\dag}T^{*}E_{T_1}, \ \ Y \succcurlyeq 0.
\label{627}
\end{align}
Under the condition $B \succcurlyeq 0$, (\ref{626}) and (\ref{627}) have solutions.  In this case,
it would be of interest to characterize the following four inequalities
\begin{align}
X \succ  Y \succcurlyeq 0, \ \  X \succcurlyeq Y \succcurlyeq 0, \ \  Y \succ  X \succcurlyeq 0, \ \   Y \succcurlyeq X \succcurlyeq 0.
\label{628}
\end{align}

As demonstrated in the previous sections,  matrix ranks and inertias and their optimizations problems are
 one of the most productive research field in matrix theory. The present author and his collaborators paid great attention in recent years for
  the development of the theory on matrix ranks and inertias, and established thousands of expansion formulas for calculating
  ranks and inertias of matrices; see, e.g., \cite{LT-jamc,LT-nla11,LT-jota,LTT,T-cal,T-ela10,T-laa10,T-laa11,T-jmaa,T-laa12,T-mcm,T-ela,
  T-mjm,T-na,T-nla,TL-lama,TL,TV} for more details. In addition, many follow-up papers by other people were also published
   on extensions and applications of  matrix rank and inertia formulas in different situations.

\noindent
Yongge Tian\\
China Economics and Management Academy\\
Central University of Finance and Economics\\
Beijing 100081, China\\
e-mail: yongge.tian@gmail.com

\end{document}